\documentclass[]{amsart}
\usepackage{color}
 \usepackage{graphicx}
\usepackage{amsmath}
\usepackage{amssymb}
\usepackage{amsfonts}
\usepackage{amsthm}

\theoremstyle{remark}
\newtheorem{remark}{Remark}[section]

\begin{document}

\title[]{Numerical scattering for the defocusing  Davey-Stewartson II equation 
for initial data with compact support}
\author{Christian Klein}
\address{Institut de Math\'ematiques de Bourgogne, UMR 5584\\
                Universit\'e de Bourgogne-Franche-Comt\'e, 9 avenue Alain Savary, 21078 Dijon
                Cedex, France\\
    E-mail Christian.Klein@u-bourgogne.fr}

\author{Nikola Stoilov}
\address{Institut de Math\'ematiques de Bourgogne, UMR 5584\\
                Universit\'e de Bourgogne-Franche-Comt\'e, 9 avenue Alain Savary, 21078 Dijon
                Cedex, France\\
    E-mail Nikola.Stoilov@u-bourgogne.fr}
\date{\today}

\begin{abstract}
In this work we present spectral algorithms for the numerical scattering for the 
defocusing
Davey-Stewartson (DS) II equation with initial data having compact 
support on a disk, i.e.,  for the solution of d-bar problems. 
Our algorithms use polar coordinates and implement a Chebychev 
spectral scheme for the radial dependence and a Fourier 
spectral method for the azimuthal dependence. The focus is placed on the 
construction of complex geometric optics (CGO) solutions which are needed in 
the scattering approach for DS. We discuss two different 
approaches: The first constructs a fundamental solution 
to the d-bar system and applies the CGO conditions on the latter. 
This is especially efficient for small values of the modulus of the 
spectral parameter $k$.  The 
second approach uses a fixed point iteration on a reformulated d-bar 
system containing the spectral parameter explicitly, a price paid to have simpler asymptotics. The approaches are illustrated for the example 
of the characteristic function of the disk and are shown to exhibit 
spectral convergence, i.e., an exponential decay of the numerical 
error with the number of collocation points. An asymptotic formula 
for large $|k|$ is given for the reflection coefficient. 
\end{abstract}
\maketitle

\section{Introduction}
Dispersive shock waves (DSW), i.e., zones of rapid modulated 
oscillations, appear in many applications in the vicinity of shocks 
whenever dispersion 
dominates dissipation. For the Korteweg-de Vries (KdV) equation, the seminal work 
by Gurevitch and Pitaevski \cite{GP} for step-like initial data provides analytic understanding of DSWs. The present work  provides the first numerical 
step in the treatment of initial data with compact support for the 
Davey-Stewartson (DS) II equation,
\begin{equation}
\begin{split}
i q_t + (q_{xx}-q_{yy}) + 2\sigma(\Phi + |q|^2)q&=0,\\
\Phi_{xx}+\Phi_{yy} +2(|q|^2)_{xx}&=0,
\end{split}
    \label{eq:DSII}
\end{equation}
an integrable 2d nonlinear Schr\"odinger equation; here $\sigma=1$ in 
the \emph{defocusing} case studied in the present paper, $\sigma=-1$ the 
in the \emph{focusing} case. The reader is referred to \cite{KS15} 
for a review on DS equations (also non-integrable cases) and a 
comprehensive list of references. We give numerical 
algorithms for the computation of the scattering transform for 
initial data with compact support on a disk for $\sigma=1$. 
We note here that, much like in the linear case of Fourier transformations, the direct and inverse problems are essentially the same. 

\subsection{D-bar systems}

An asymptotic description of DSWs is in 
general only possible for small amplitudes via multiscales 
approximations. For completely integrable equations in one spatial 
dimension, steepest descent techniques for Riemann-Hilbert problems  
(RHPs) 
allow results for large amplitudes as well. A complete description of 
DSWs including the appearence of Painlev\'e transcendents exists so 
far only for the Korteweg-de Vries equation, see \cite{GK} for a 
review. 
For two-dimensional integrable equations d-bar problems often take 
the role RHPs play the 1d case, see \cite{BC}. For such 
problems no steepest descent techniques are known yet though partial 
progress 
has been made recently in \cite{AKMM} for the defocusing 
Davey-Stewartson II equation.  In any case so far DSWs for 
two-dimensional integrable equations have been mainly studied 
numerically, see \cite{KSM,KR13,KR14}. 

The direct scattering transformation for the  DS II equations is 
given by a system of d-bar equations,
\begin{equation}
\begin{array}{l}
    \bar{\partial}\psi_{1}=\frac{1}{2}q\psi_{2},\\
    ~\\
    \partial\psi_{2}=\frac{1}{2}\sigma\bar{q}\psi_{1},
    \label{dbarpsi}
\end{array}
\end{equation}
where the scalar functions $\psi_{1}$ and $\psi_{2}$ satisfy the 
 \emph{complex geometric 
optics} (CGO) asymptotic conditions
\begin{equation}
\begin{array}{l}
    \lim_{|z|\to\infty}\psi_{1}\mathrm{e}^{-kz}=1,\\
    ~\\
    \lim_{|z|\to\infty}\psi_{2}\mathrm{e}^{-\bar{k}\bar{z}}=0;
    \label{dbarpsiasym}
\end{array}
\end{equation}
here  $q=q(x,y,t)$ is a complex-valued field, the \emph{spectral 
parameter} $k\in\mathbb{C}$ is independent of $z=x+\mathrm{i} y$, and
\begin{equation*}
\partial:=\frac{1}{2}\left(\frac{\partial}{\partial x}-\mathrm{i}\frac{\partial}{\partial y}\right)\quad\text{and}\quad
\bar{\partial}:=\frac{1}{2}\left(\frac{\partial}{\partial 
x}+\mathrm{i}\frac{\partial}{\partial y}\right).
\end{equation*}
The \emph{reflection coefficient} $R=R(k)$ (which can be seen as a 
nonlinear analogue to a Fourier transform) is defined in terms of 
$\psi_2(z;k)$\footnote{Note that the notation $\psi_{2}(z;k)$ does 
not imply that the function is holomorphic in either $z$ or $k$.} as follows:
\begin{equation}
\mathrm{e}^{-kz}\overline{\psi_2(z;k)}=\frac{1}{2}R(k) z^{-1}+ O(|z|^{-2}),\quad |z|\to\infty.
\label{eq:r-def}
\end{equation}
The existence and uniqueness of CGO solutions to system 
(\ref{dbarpsi}) with $\sigma=1$ was studied in \cite{Su1,Su2}) for Schwartz class 
potentials and in \cite{BU,Brown,Perry2012,NRT} for more general 
potentials. Note that the understanding is much less complete in the 
focusing case $\sigma=-1$ since the system (\ref{dbarpsi}) no longer 
has generically a unique solution for large classes of potentials $q$ 
for all $k\in\mathbb{C}$. There can 
be so-called exceptional points (special values of the spectral parameter 
$k$) where the system is not uniquely solvable. Note that the codes
presented in this paper can also be applied to the focusing case in 
the absence of such exceptional points. However, since it is not clear when 
these points appear and since numerical problems are expected in the 
vicinity of exceptional points, we concentrate on the defocusing case 
here.  

Systems of the form (\ref{dbarpsi}) also appear in electrical impedance tomography 
(EIT) in 2d, the reconstruction of 
the conductivity in a given domain from measurements of the electrical 
current through its boundary, induced by an applied voltage, i.e., from 
the Dirichlet-to-Neumann map. This problem was first posed by Calderon
\cite{calderon} and bears his name. For a comprehensive review of the 
mathematical aspects and advances see \cite{Uhl}. They also appear
 in the context of 2d orthogonal polynomials, and of 
Normal Matrix Models in Random Matrix Theory, see e.g.\ \cite{KM}.

As $q$ in (\ref{eq:DSII}) evolves in time $t$ according to \eqref{eq:DSII}, the reflection coefficient evolves by a trivial phase factor:
\begin{equation}
R(k;t)=R(k,0)e^{4i t\Re(k^2)}.
\end{equation}
The inverse scattering transform for DS II is then given by 
(\ref{dbarpsiasym}) after replacing $q$ by $R$ and vice versa, the derivatives with 
respect to $z$ by the corresponding derivatives with respect to $k$, 
and asymptotic conditions for $k\to\infty$ instead of $z\to\infty$. 

\subsection{Numerical appoaches}

Numerical approaches to d-bar systems have so far mainly taken the 
following path: the inverse of the d-bar operator is known to be 
given by the solid Cauchy transform 
\begin{equation}
\overline{\partial}^{-1}F(x,y):=-\frac{1}{\pi}\int_{}^{}\int_{\mathbb{R}^2}\frac{F(x',y')\,dx'dy'}{(x'-x)+\mathrm{i} (y'-y)},
\label{eq:solid-Cauchy}
\end{equation}
a weakly singular integral. This integral can be either computed 
via finite difference, finite element 
discretizations, or via a Fourier approach, see \cite{MuSi2012} for a 
review of techniques. 

The original approach \cite{KnMuSi2004}, see also \cite{APRS},  uses a Fourier method for the computation of the integral, 
i.e., a 2d discrete Fourier transform. A
rather bold regularization of the singular integrand was implemented by replacing the 
singularity simply by a finite value. This leads to a finite, but 
discontinuous integrand. Though discrete Fourier transforms are 
spectral methods which are known to show spectral convergence, i.e. 
an exponential decrease of the numerical error with the resolution, for the 
approximation of smooth functions, they are of first order 
method if a non-continuous function approximated. This means that the 
numerical error decreases only linearly with $1/N$ where $N$ is the 
number of Fourier modes. This  first order convergence was proven in 
\cite{KnMuSi2004}.

The first Fourier approach with spectral convergence (i.e., 
exponential decrease of the numerical error with the number $N$ of 
Fourier modes)
was presented in \cite{KM} for potentials in the 
Schwartz class of rapidly decreasing smooth potentials $q$. The 
approach uses an analytic (up to numerical precision) regularization 
of the integrand.  This is further developed in \cite{KMS}. In 
practice this means that accuracies of the order of machine precision 
can be reached on low-cost computers which are out of reach of a 
first order method. 

In various applications,  treating potentials with compact support is crucial. 
The focus of this paper will be on such potentials. 
The method of \cite{KM} cannot be generalized to potentials with 
compact support without severe loss of spectral accuracy since Gibbs phenomena appear at the discontinuity, see for example \cite{APRS} where a Fourier solver for the Beltrami equation is applied to various step like potentials. 
Therefore, here  we present a completely different 
approach starting from a formulation of the system (\ref{dbarpsi}) in 
polar coordinates. The system is discretized with a Chebychev 
spectral method in the radial coordinate $r$ and a Fourier spectral 
method in the angular coordinate $\phi$. This means that the functions 
$\psi_{1,2}$ in (\ref{dbarpsi}) will be approximated via truncated 
series of Chebychev polynomials in $r$ and truncated Fourier series in $\phi$. 
Fast, efficient algorithms exist for both these approaches, making them a prime choice. In addition they are known to show an 
exponential decrease of the numerical error with the number of 
Chebychev polynomials and Fourier modes for smooth functions, i.e., 
here for smooth potentials $q$ on the disk. We show at examples that 
no Gibbs phenomena can be observed in our approach since we only 
approximate smooth functions (the functions are smooth on the 
considered domain, discontinuities appear only on the boundary). 
The numerical errors we discuss are 
thus always global, even at the discontinuity of the potential.

This discretization will be used in two different ways: First we 
construct a fundamental system of solutions to the finite dimensional 
system of ordinary differential equations (ODEs) with which 
(\ref{dbarpsi}) is approximated after discretization in $\phi$. Note 
that the discontinuity of $q$ at the boundary of the disk does not 
affect the method since the latter is a domain boundary. Smoothness 
is only required inside the computational domain. 
The approach 
can be applied for arbitrary $k$ since the system (\ref{dbarpsi}) is 
independent of the spectral parameter $k$. The CGO conditions 
(\ref{dbarpsiasym}) are then implemented for given values of $k$ by 
imposing them on this fundamental solution. It is shown that this 
works well for values of $|k|\lesssim 1$. 

For larger $|k|$, cancellation errors will play a role in this 
approach. To address 
this, we introduce functions $\Phi_{1}=\mathrm{e}^{-kz}\psi_{1}$ and 
$\Phi_{2}=\mathrm{e}^{-\bar{k}\bar{z}}\psi_{2}$ which satisfy with 
(\ref{dbarpsiasym}) the simpler asymptotic conditions 
\begin{equation}
    \lim_{|z|\to\infty}\Phi_{1}=1,\quad \lim_{|z|\to\infty}\Phi_{2}=0.
    \label{Phisasym}
\end{equation}
In other words, $\Phi_{1}$ and $\Phi_{2}$ are bounded functions at 
infinity. In terms of $\Phi_{1}$ and $\Phi_{2}$ the d-bar system 
(\ref{dbarpsi}) leads to
\begin{equation*}
\begin{array}{l}
    \bar{\partial}\Phi_{1}=\frac{1}{2}q\mathrm{e}^{\bar{k}\bar{z}-kz}\Phi_{2},\\
    ~\\
    \partial\Phi_{2}=\frac{1}{2}\bar{q}\mathrm{e}^{kz-\bar{k}\bar{z}}\Phi_{1},
\end{array}
\end{equation*}
for $\Phi_{1}$, $\Phi_{2}$ which contains the spectral parameter explicitly. The 
latter is solved with the same Chebychev and Fourier discretization 
as before, but this time with a fixed point iteration. This iteration is 
shown to converge rapidly in examples for all $k$ which allows to 
treat efficiently large values of $|k|$ and thus rapidly oscillating solutions. 

The paper is organized as follows: In Section 2 we give analytical justification of the equations we are going to treat numerically using polar coordinates and study the example of 
$q$ being the characteristic function of the disk. In Section 3 we 
summarize the numerical schemes employed to discretize the 
differential operators and construct a fundamental solution to 
(\ref{dbarpsi}). This fundamental solution is used in Section 4 to 
construct CGO solutions. In Section 5 we apply a fixed point iteration to 
solve system (\ref{dbarpolphi}) which is equivalent to 
(\ref{dbarpsi}). In Section 6 we add some concluding remarks.

\section{D-bar system with compactly supported potentials on a disk}
In this section we formulate the defining equations (\ref{dbarpsi}) 
and (\ref{dbarpsiasym}) for the CGO 
solutions in polar coordinates and use Fourier series in the 
azimuthal variable for the solutions. As an example we consider the 
case of constant potential on the disk.

\subsection{Polar coordinates}

We write $z=r\mathrm{e}^{\mathrm{i\phi}}$ in which (\ref{dbarpsi}) reads 
\begin{equation}
\begin{array}{l}
    \mathrm{e}^{\mathrm{i}\phi}\left(\partial_{r}+\frac{\mathrm{i}}{r}\partial_{\phi}\right)\psi_{1}  
    =q(r,\phi)\psi_{2},
\\
~\\
     \mathrm{e}^{-\mathrm{i}\phi}\left(\partial_{r}-\frac{\mathrm{i}}{r}\partial_{\phi}\right)\psi_{2}  
    =\bar{q}(r,\phi)\psi_{1}.
    \label{dbarpol}
\end{array}
\end{equation}We are looking for solutions in the form of a Fourier series,
\begin{equation}
    \psi_{1}=\sum_{n\in\mathbb{Z}}^{}a_{n}(r)\mathrm{e}^{\mathrm{i}n\phi},\quad 
     \psi_{2}=\sum_{n\in\mathbb{Z}}^{}b_{n}(r)\mathrm{e}^{\mathrm{i}n\phi}.
    \label{diskasym}
\end{equation}
The functions $a_{n}$ and $b_{n}$ are determined by (\ref{dbarpol}) 
for given $q(r,\phi)$
together with a regularity condition at $r=0$ up to a certain number 
of constants $\alpha_{n}$ and $\beta_{n}$ corresponding to general 
holomorphic respectively antiholomorphic functions arising in the determination 
of $a_{n}$ and $b_{n}$ respectively: if $a_{n}$ and $b_{n}$ are 
solutions to (\ref{dbarpol}) regular for $r\to0$, so are
\begin{equation}
    a_{n}(r)+\alpha_{n}r^{n},\quad n=0,1,2,\ldots,\quad 
    b_{n}(r)+\beta_{n}r^{n},\quad n=0,-1,-2,\ldots
    \label{ab}
\end{equation}
The constants $\alpha_{n}$, $\beta_{n}$ will be uniquely determined by the 
asymptotic conditions (\ref{dbarpsiasym}).

For $r>1$, where $q\equiv0$, the function $\psi_{1}$ is a holomorphic function, whereas $\psi_{2}$ is
antiholomorphic. Continuity of the functions at the disk thus implies 
for $r>1$
\begin{equation*}
    \psi_{1}=\sum_{n\in\mathbb{Z}}^{}a_{n}(1)r^{n}\mathrm{e}^{\mathrm{i}n\phi},\quad 
     \psi_{2}=\sum_{n\in\mathbb{Z}}^{}b_{n}(1)r^{-n}\mathrm{e}^{\mathrm{i}n\phi}.
\end{equation*}
The first of the asymptotic conditions (\ref{dbarpsiasym}) is 
equivalent to  $c_{0}=1$ and $c_{n}=0$ for $n>0$ for $r\to\infty$, where 
\begin{equation}
    \sum_{n\in\mathbb{Z}}^{}c_{n}r^{n}\mathrm{e}^{\mathrm{i}n\phi}:=
    \sum_{n=0}^{\infty}\sum_{m=0}^{n}\frac{(-k)^{m}}{m!}\left(a_{n-m}(1)r^{n}\mathrm{e}^{\mathrm{i}n\phi}
    +a_{m-n-1}(1)r^{2m-n-1}\mathrm{e}^{\mathrm{i}(2m-n-1)\phi}\right)
    \label{diskasym1}.
\end{equation}
Since the conditions are imposed for $r\to\infty$, negative powers of 
$r$ in (\ref{diskasym1}) can be neglected. This implies that there 
will be only conditions for the positive frequencies $c_{n}$, $n\geq0$, which 
can be written in the form, 
\begin{equation}
    c_{n}=\sum_{m=0}^{\infty}a_{n-m}(1)\frac{(-k)^{m}}{m!}
    \label{cn}.
\end{equation}
Note that these conditions --- though obtained in the limit 
$r\to\infty$ --- lead here to conditions at the rim of the disk  $r=1$ because of the 
holomorphicity of $\psi_{1}$ for $r>1$ and the continuity of the 
function at the rim of the disk. In other words, the solution in the 
exterior of the disk follows from the holomorphicity of $\psi_{1}$ 
there as well as the asymptotic condition and the continuity at the 
disk. It is not necessary to solve a PDE there in contrast to the 
disk. The convolutions in 
(\ref{cn}) can be computed in  a standard way via Fourier series,
\begin{equation}
   c_{n} =\frac{1}{2\pi}\int_{0}^{2\pi}\mathrm{d}\phi
   \exp\left(-k\mathrm{e}^{\mathrm{i\phi}}\right)
   \sum_{m\in\mathbb{Z}}^{}a_{m}(1)\mathrm{e}^{\mathrm{i(m-n)\phi}}
    \label{cnfourier}.
\end{equation}
In the numerical approach we  present in Section 4, we will use the asymptotic conditions in the form (\ref{cnfourier}).

Similarly we get for the second condition in (\ref{dbarpsiasym}) that 
$d_{n}=0$ for $n=0,-1,-2,\ldots$, where
\begin{equation*}
    \sum_{n\in\mathbb{Z}}^{}d_{n}r^{n}\mathrm{e}^{\mathrm{i}n\phi}:=
    \sum_{n=0}^{\infty}\sum_{m=0}^{n}\frac{(-\bar{k})^{m}}{m!}
    \left(b_{n-m+1}(1)r^{2m-n-1}\mathrm{e}^{\mathrm{i}(n+1-2m)\phi}
    +    
    b_{m-n}(1)r^{n}\mathrm{e}^{-\mathrm{i}n\phi}\right)
\end{equation*}
i.e., 
\begin{equation*}
    d_{n} = \sum_{m=0}^{\infty}b_{m-n}(1)\frac{(-\bar{k})^{m}}{m!}
\end{equation*}
This is equivalent to 
\begin{equation}
   d_{n} 
   =\frac{1}{2\pi}\int_{0}^{2\pi}\mathrm{d}\phi\exp\left(-\bar{k}\mathrm{e}^{\mathrm{-i\phi}}\right)
   \sum_{m\in\mathbb{Z}}^{}b_{m}(1)\mathrm{e}^{\mathrm{i(n-m)\phi}}
    \label{dnfourier}.
\end{equation}

In this approach the reflection 
coefficient (\ref{eq:r-def}) is simply given by the coefficient 
$d_{1}$,
\begin{equation}
    R = 2\bar{d}_{1}.
    \label{Rd}
\end{equation}

System (\ref{dbarpsi}) has the advantage that it is independent of 
the spectral parameter $k$. A disadvantage from a numerical point of 
view are the asymptotic conditions (\ref{dbarpsiasym}) which imply 
that the functions $\psi_{1}$ and $\psi_{2}$ have essential 
singularities at infinity. An alternative way to treat 
(\ref{dbarpsi}) is thus to introduce the functions
\begin{equation}
    \Phi_{1}=\mathrm{e}^{-kz}\psi_{1},\quad 
    \Phi_{2}=\mathrm{e}^{-\bar{k}\bar{z}}\psi_{2}
    \label{Phi},
\end{equation}
which satisfy the asymptotic conditions (\ref{Phisasym}).
In the case of a potential $q$ with compact support, 
$\Phi_{1}$ is holomorphic in the complement of the support of $q$ and 
$\Phi_{2}$ is antiholomorphic there. System (\ref{dbarpol}) for the 
functions (\ref{Phi}) reads 
\begin{equation}
\begin{array}{l}
    \mathrm{e}^{\mathrm{i}\phi}\left(\partial_{r}+\frac{\mathrm{i}}{r}\partial_{\phi}\right)\Phi_{1}  
    =q(r,\phi)\mathrm{e}^{\bar{k}\bar{z}-kz}\Phi_{2},
    \\~\\
     \mathrm{e}^{-\mathrm{i}\phi}\left(\partial_{r}-\frac{\mathrm{i}}{r}\partial_{\phi}\right)\Phi_{2} 
    =\bar{q}(r,\phi)\mathrm{e}^{kz-\bar{k}\bar{z}}\Phi_{1}.
    \label{dbarpolphi}
\end{array}
\end{equation}
This system contains, in contrast to (\ref{dbarpol}), not just $k$ 
but (for $|k|$ large) rapidly oscillating functions. 
On the other hand, the asymptotic conditions (\ref{Phisasym}) are 
considerably simpler than the conditions (\ref{dbarpsiasym}). 
Therefore we will use (\ref{dbarpol}) for values of $|k|\lesssim 1$, 
and (\ref{dbarpolphi}) for larger values of $|k|$.

\subsection{Example}\label{example}
In order to illustrate the above approach and to have a concrete 
example to test at least parts of the code, we consider the case of a 
constant potential at the disk, 
\begin{equation*}
       q = 
    \begin{cases}
        1, & r\leq 1 \\
        0, & r>1
    \end{cases}. 
\end{equation*}

Differentiating the first equation of (\ref{dbarpsi}) with respect to 
$\partial$ and eliminating $\psi_{2}$ via the second, we get in polar coordinates for $r\leq 1$
\begin{equation*}
    \partial_{rr}\psi_{1}+\frac{1}{r}\partial_{r}\psi_{1}
    +\frac{1}{r^{2}}\partial_{\phi\phi}\psi_{1}=\psi_{1}.   
\end{equation*}
Note that the same equation holds for $\psi_{2}$. With 
(\ref{diskasym}), we get that 
the $a_{n}$ satisfy the modified Bessel equation
\begin{equation*}
    a_{n}''+\frac{1}{r}a_{n}'-\frac{r^{2}+n^{2}}{r^{2}}a_{n}=0,
\end{equation*}
i.e., the  solution regular on the disk is given by 
$a_{n}=\alpha_{n}I_{n}(r)$, where the $\alpha_{n}=const$ and where 
the $I_{n}(r)$ are the modified Bessel functions \cite{AS}. Thus we 
have for 
$r\leq 1$
\begin{equation}
    \psi_{1}=\sum_{n\in\mathbb{Z}}^{}\alpha_{n}I_{n}(r)\mathrm{e}^{\mathrm{i}n\phi}
    \label{psi11}.
\end{equation}

For $r\leq 1$, the function $\psi_{2}$ can be computed from the first 
relation of (\ref{dbarpsi}).
Because of the identity $I_{n}'(r)-(n/r)I_{n}=I_{n+1}$, it can be written in the form
\begin{equation}
    \psi_{2}=        
   \sum_{n\in\mathbb{Z}}^{}\alpha_{n}\mathrm{e}^{\mathrm{i}(n+1)\phi}I_{n+1}(r).
    \label{psi22}
\end{equation}

For $r>1$, because of the continuity of the potentials at $r=1$  we get
\begin{equation}
 \psi_{1}=\sum_{n\in\mathbb{Z}}^{}\alpha_{n}I_{n}(1)z^{n},\quad 
      \psi_{2}=        
   \sum_{n\in\mathbb{Z}}^{}\frac{\alpha_{n}}{\bar{z}^{n+1}}I_{n+1}(1).
    \label{psi12}
\end{equation}

The constants $\alpha_{n}$, $n\in\mathbb{Z}$ are determined by the asymptotic 
condition (\ref{dbarpsiasym}). We only found an explicit solution for 
$k=0$, 
$\alpha_{0}=1/I_{0}(1)$ whereas $\alpha_{n}=0$ for $n\neq 0$. This 
implies
\begin{equation*}
       \psi_{1} = 
    \begin{cases}
        I_{0}(r)/I_{0}(1), & r\leq 1 \\
        1, & r>1
    \end{cases}, 
\end{equation*}
and
\begin{equation*}
       \psi_{2} = 
    \begin{cases}
        I_{1}(r)/I_{0}(1)\mathrm{e}^{\mathrm{i}\phi}, & r\leq 1 \\
        I_{1}(1)/I_{0}(1)/\bar{z}, & r>1
    \end{cases}, 
\end{equation*}

\section{Numerical construction of a fundamental solution to the 
d-bar system}
In this section we present a spectral approach to the system 
(\ref{dbarpol}) based on a discrete Fourier approach in $\phi$ and a 
Chebychev spectral method in $r$. We construct solutions 
$\psi_{1,2}^{(j)}$ to the system characterized by  boundary 
conditions on the Fourier coefficients (\ref{diskasym}) at the rim of 
the disk, in a way that a basis of solutions is obtained. In the next section this 
system of fundamental solutions is subjected to the asymptotic conditions (\ref{cnfourier}) and 
(\ref{dnfourier}). The solutions are then tested for the example of 
subsection \ref{example}.

\subsection{Spectral approach}
The periodicity in the coordinate $\phi$ in (\ref{dbarpol}) suggests 
to use Fourier techniques for this variable. This means that the 
series in (\ref{diskasym}) are approximated via discrete Fourier 
transforms which implies the use of a discrete variable $\phi$ 
sampled on the $N_{\phi}$ collocation points 
\begin{equation*}
    \phi_{j}= \frac{2\pi j}{N_{\phi}},\quad j=0,1,\ldots,N_{\phi}-1
\end{equation*}
where $N_{\phi}$ is an even positive integer. Note that the discrete 
Fourier series is not only periodic in $\phi$, but also in the dual 
variable $n$. Since the d-bar system (\ref{dbarpol}) contains factors 
$\mathrm{e}^{\pm\mathrm{i}\phi}$, we use an asymmetric definition of 
the Fourier sums approximating (\ref{diskasym}),
\begin{equation}
    \psi_{1}\approx \sum_{n=-N_{\phi}/2}^{N_{\phi}/2-1}a_{n}(r)\mathrm{e}^{\mathrm{i}n\phi},\quad 
     \psi_{2}\approx\sum_{n=-N_{\phi}/2+1}^{N_{\phi}/2}b_{n}(r)\mathrm{e}^{\mathrm{i}n\phi},
    \label{diskfft}
\end{equation}
in order to be able to deal with the same powers of 
$\mathrm{e}^{\mathrm{i}\phi}$ in both equations (we use the same 
symbols as in (\ref{dbarpol}) to avoid cluttered notation). The 
discrete Fourier transforms can be computed conveniently with a 
\emph{Fast Fourier Transform} (FFT). 

With (\ref{diskfft}), the partial differential equation (PDE) system 
(\ref{dbarpol}) is approximated via a finite system of ordinary 
differential equations (ODE) in $r$ for the $2N_{\phi}$ functions 
$a_{n}(r)$, $n=-N_{\phi}/2,\ldots,N_{\phi}/2-1$, and $b_{n}(r)$, 
$n=-N_{\phi}/2+1,\ldots,N_{\phi}/2$. To solve this system, we 
approximate these functions via a sum of Chebychev polynomials 
$T_{m}(l):=\cos (m\arccos x)$, $l\in[-1,1]$, i.e., 
\begin{equation}
\begin{array}{l}
    a_{n}(r)\approx \sum_{m=0}^{N_{r}}a_{nm} T_{m}(l),\quad 
    n=-\frac{N_{\phi}}{2},\ldots,\frac{N_{\phi}}{2}-1,\\~\\
    b_{n}(r)\approx \sum_{m=0}^{N_{r}}b_{nm} T_{m}(l),\quad 
    n=-\frac{N_{\phi}}{2}+1,\ldots,\frac{N_{\phi}}{2}
    \label{cheb},
\end{array}
\end{equation}
where $r=(1+l)/2$. 
The constants $a_{nm}$ are determined via a collocation method: on 
the points 
\begin{equation*}
    l_{j}=\cos\left(\frac{j\pi}{N_{r}}\right), \quad j=0,\ldots,N_{r}
\end{equation*}
the equations (\ref{cheb}) 
are imposed as equalities, i.e., for fixed $n$ we have
\begin{equation*}
    a_{n}(r(l_{j})) = \sum_{m=0}^{N_{r}}a_{nm}T_{m}(l_{j}),\quad 
    j=0,\ldots,N_{r}
\end{equation*}
which uniquely determines the $a_{nm}$. In the same way the $b_{nm}$ 
are fixed. Because of the definition of the Chebychev polynomials, 
one has
$T_{m}(l_{j})=\cos (mj\pi/N_{r})$. Thus the coefficients $a_{nm}$ and 
$b_{nm}$ can be obtained via a \emph{Fast Cosine Transformation} 
(FCT) which is related to the FFT, see e.g., the discussion in 
\cite{trefethen} and references therein. 

The existence of fast algorithms, however, is not the only 
advantage of the Chebychev spectral method and the discrete 
Fourier transform. Both are so-called spectral methods which means 
that the numerical error in approximating smooth functions decreases 
faster than any power of $1/N_{r}$ and $1/N_{\phi}$, in practice it 
decreases exponentially. This is due to an analogue for discrete 
Fourier transforms of the well known theorem that 
the Fourier transform of a smooth, rapidly decreasing function is 
rapidly decreasing, see the discussion in \cite{trefethen}. 

To approximate derivatives via the ansatz (\ref{cheb}), one uses 
$T_{0}'(l)=0$, $T_{1}'(l)=1$ and for $n\geq1$ the 
identity 
\begin{equation*}
  \frac{T'_{n+1}(l)}{n+1} - \frac{T'_{n-1}(l)}{n-1} = 2 T_n(l)
\end{equation*}
which implies that the derivative of the $a_{n}(r)$ is approximated via the 
action of a \emph{differentiation matrix} $D$ on the Chebychev coefficients 
$a_{nm}$ 
\begin{equation*}
    a_{n}'(r)\approx 
    \sum_{m,\alpha=0}^{N_{r}}D_{m\alpha}a_{n\alpha}T_{m}(l)
\end{equation*}
The differentiation matrix is upper triangular and for even $N_{r}$ of the form 
$$D = 
\begin{pmatrix}
    0 & 1 & 0 & 3 & 0 & 5 & \hdots & N_{r}-1&0 \\
    ¥ & ¥ & 4 & 0 & 8 & \hdots & ¥ & ¥ & 2*N_{r}\\
    ¥ & ¥ & ¥ & 6 & 0 & 10 & \hdots & 2*(N_{r-1}) & 0\\
    ¥ & ¥ & ¥ & ¥ & 8 & 0 & 12 & \hdots & 2*N_{r} \\
    ¥ & ¥ & ¥ & ¥ & ¥ & \ddots & ¥ & ¥&0 \\
    ¥ & ¥ & ¥ & ¥ & ¥ & ¥ & ¥ & ¥ &\vdots\\
    ¥ & ¥ & ¥ & ¥ & ¥ & ¥ & ¥ & ¥ &0\\
    ¥ & ¥ & ¥ & ¥ & ¥ & ¥ & ¥ & ¥& 2*N_{r}
\end{pmatrix}.
$$

In a similar way one can divide in the space of Chebychev 
coefficients by $r$ which is normally numerically a very delicate 
operation if $r$ can vanish on the considered interval as it does 
here. Using the identity
\begin{equation}
    T_{n+1}(l)+T_{n-1}(l)=2lT_{n}(x), \quad n=1,2,\ldots
    \label{chebrec},
\end{equation}
we can divide in coefficient space by $l\pm1$. 
For given Chebyshev coefficients $a_{nm}$ we define the coefficients 
$\tilde{a}_{nm}$ 
via 
$\sum_{n=0}^{\infty}a_{nm}T_{m}(l)=:\sum_{n=0}^{\infty}(l\pm1)\tilde{a}_{nm}T_{m}(l)$. 
This implies the action of a matrix $R$ in coefficient space (if 
$a_{n}(r)/r$ is bounded for $r\to0$),
\begin{equation}
    \frac{a_{n}(r)}{r}\approx  
    \sum_{m,\alpha=0}^{N_{r}}R_{m\alpha}a_{n\alpha}T_{m}(l)
    \label{chebdiv}.
\end{equation}
The matrix $R^{-1}$ has for even $N_{r}$ the form
$$R^{-1} = 
\begin{pmatrix}
    1 & 1/2 & ¥ & ¥ & ¥ \\
    1 & 1 & 1/2 & ¥ & ¥  \\
    ¥ & 1/2 & \ddots & \ddots & \\
    ¥ & ¥ & \ddots & ¥ & ¥ \\
     ¥ & ¥ & ¥ & ¥ & ¥ & 1/2 \\
     ¥ & ¥ & ¥ & ¥ & 1/2 & 1
\end{pmatrix};
$$
(the matrix $R$ is computed by inverting this matrix).

Denoting the combined action of FFT and FCT on a function 
$\psi(r,\phi)$ by $\mathbb{F}$, we thus 
approximate the system (\ref{dbarpol}) via ($m=0,\ldots,N_{r}$)
\begin{align}
    \sum_{\alpha=0}^{N_{r}}\left(D-nR\right)_{m\alpha}a_{n\alpha}&=\mathbb{F}\left(q(r,\phi) \mathrm{e}^{-\mathrm{i}\phi}
    \sum_{\alpha=0}^{N_{r}}\sum_{\beta=-N_{\phi}/2}^{N_{\phi}/2-1}b_{\beta\alpha}T_{\alpha}(l)\mathrm{e}^{\mathrm{i}\beta\phi}\right)_{nm},\quad n=-N_{\phi}/2,\ldots,N_{\phi}/2-1
    \nonumber\\
        \sum_{\alpha=0}^{N_{r}}\left(D+nR\right)_{m\alpha}b_{n\alpha}&=\mathbb{F}\left(\bar{q}(r,\phi) \mathrm{e}^{\mathrm{i}\phi}
    \sum_{\alpha=0}^{N_{r}}\sum_{\beta=-N_{\phi}/2}^{N_{\phi}/2-1}a_{\beta\alpha}T_{\alpha}(l)\mathrm{e}^{\mathrm{i}\beta\phi}\right)_{nm},\quad n=-N_{\phi}/2+1,\ldots,N_{\phi}/2.
    \label{dbarsysnum}
\end{align}
The action of the functions $q(r,\phi)\mathrm{e}^{-\mathrm{i}\phi}$ and its 
conjugate in (\ref{dbarsysnum}) is computed as a convolution in 
coefficient space via (\ref{chebrec}) and by using the periodicity of 
the FFT in Fourier space. Note that the computation of convolutions 
or products in a spectral method can lead to so called 
\emph{aliasing} errors, see the discussion in \cite{trefethen}, which 
are due to the fact that only a finite number of terms are considered 
in the Chebychev and Fourier sums, whereas for instance 
relation (\ref{chebrec}) implies the use of the whole series in order 
to get equality. 

Denoting the vector built from the  coefficients 
$a_{nm}$ by $A$ and the vector built by the coefficients $b_{nm}$ by 
$B$, equation (\ref{dbarsysnum}) can be formally written as 
\begin{equation}
    \mathcal{O}
    \begin{pmatrix}
        A \\
        B
    \end{pmatrix}=0
    \label{mat},
\end{equation}
where $\mathcal{O}$ is the 
$2(N_{r}+1)N_{\phi}\times2(N_{r}+1)N_{\phi}$ matrix 
built from the matrices 
$D$, $R$ and the convolutions of the coefficients of 
$q\mathrm{e}^{-\mathrm{i}\phi}$ and its conjugate. 

\subsection{Fundamental solution}
The matrix $\mathcal{O}$ in (\ref{mat}) has an $N_{\phi}$-dimensional 
kernel corresponding to the homogeneous solutions characterized by 
the constants $\alpha_{n}$, $\beta_{n}$ in (\ref{ab}). Note that 
regularity of the solutions at $r=0$ does not need to be imposed 
since we handle the terms proportional to $1/r$ in (\ref{dbarpol}) in the space of 
coefficients via the matrix $R$ in (\ref{chebdiv}). 

To obtain the general solution to the system (\ref{dbarsysnum}), we 
introduce functions $\psi^{(j)}_{1}$ and $\psi^{(j)}_{2}$, 
$j=1,2,\ldots,N_{\phi}$ such that the general solution to the system 
can be written in the form 
\begin{equation}
    \psi_{1}=\sum_{j=1}^{N_{\phi}}\gamma_{j}\psi^{(j)}_{1},\quad 
    \psi_{2}=\sum_{j=1}^{N_{\phi}}\gamma_{j}\psi^{(j)}_{2}
    \label{general},
\end{equation}
where $\gamma_{j}=const$ for $j=1,\ldots,N_{\phi}$. We define the 
solutions $\psi^{(j)}_{1}$ and $\psi^{(j)}_{2}$ to (\ref{dbarsysnum}) 
in the following way via boundary conditions at the rim of the disk: \\
for $j=1,\ldots,N_{\phi}/2$
\begin{equation}
    a_{n}(1)=\delta_{(j-1)n},\quad n=-N_{\phi}/2,\ldots, 
    N_{\phi}/2-1, \quad b_{n}(1)=0, \quad n=-N_{\phi}/2+1,\ldots, 
    N_{\phi}/2
    \label{abcond1},
\end{equation}
and for $j=N_{\phi}/2+1,\ldots,N_{\phi}$
\begin{equation}
    a_{n}(1)=0,\quad n=-N_{\phi}/2,\ldots, 
    N_{\phi}/2-1, \quad b_{n}(1)=\delta_{(j-N_{\phi})n}, \quad n=-N_{\phi}/2+1,\ldots, 
    N_{\phi}/2
    \label{abcond2}.
\end{equation}
These conditions uniquely determine the functions $\psi^{(j)}_{1}$ 
and $\psi^{(j)}_{2}$ and ensure that they form a basis of the regular 
solutions to (\ref{dbarsysnum}). 

The conditions (\ref{abcond1}) and (\ref{abcond2}) are implemented 
in the numerical approach via Lanczos' $\tau$-method \cite{tau}. The 
idea is that these conditions replace certain equations in 
(\ref{mat}) which leads to a new system of the form
\begin{equation}
    \tilde{\mathcal{O}}
    \begin{pmatrix}
        A \\
        B
    \end{pmatrix}=S
    \label{mat2};
\end{equation}
here $\tilde{O}$ is the matrix $\mathcal{O}$ where the rows 
corresponding to the  Fourier index $j$ appearing in (\ref{abcond1}) 
and (\ref{abcond2}) and the Chebyshev index $m=N_{r}+1$ (the terms 
corresponding to the highest Chebyshev polynomial are thus in these 
cases 
neglected) are replaced by the left hand sides of (\ref{abcond1}) and 
(\ref{abcond2}). The right hand side 
$S$ of (\ref{mat2}) is a $2(N_{r}+1)N_{\phi}\times N_{\phi}$ matrix which is 
identical to zero except for a single value in each column corresponding 
to the Fourier index $j$ and the Chebyshev index $N_{r}+1$. 

Solving (\ref{mat2}) thus provides a basis of solutions to 
(\ref{dbarsysnum}). The $\tau$ method ensures that the boundary 
conditions at the disk are satisfied with the same spectral accuracy 
as the solved differential equation, see the discussion in 
\cite{trefethen}. 

\subsection{Test}
\begin{figure}[htb!]
  \includegraphics[width=0.49\textwidth]{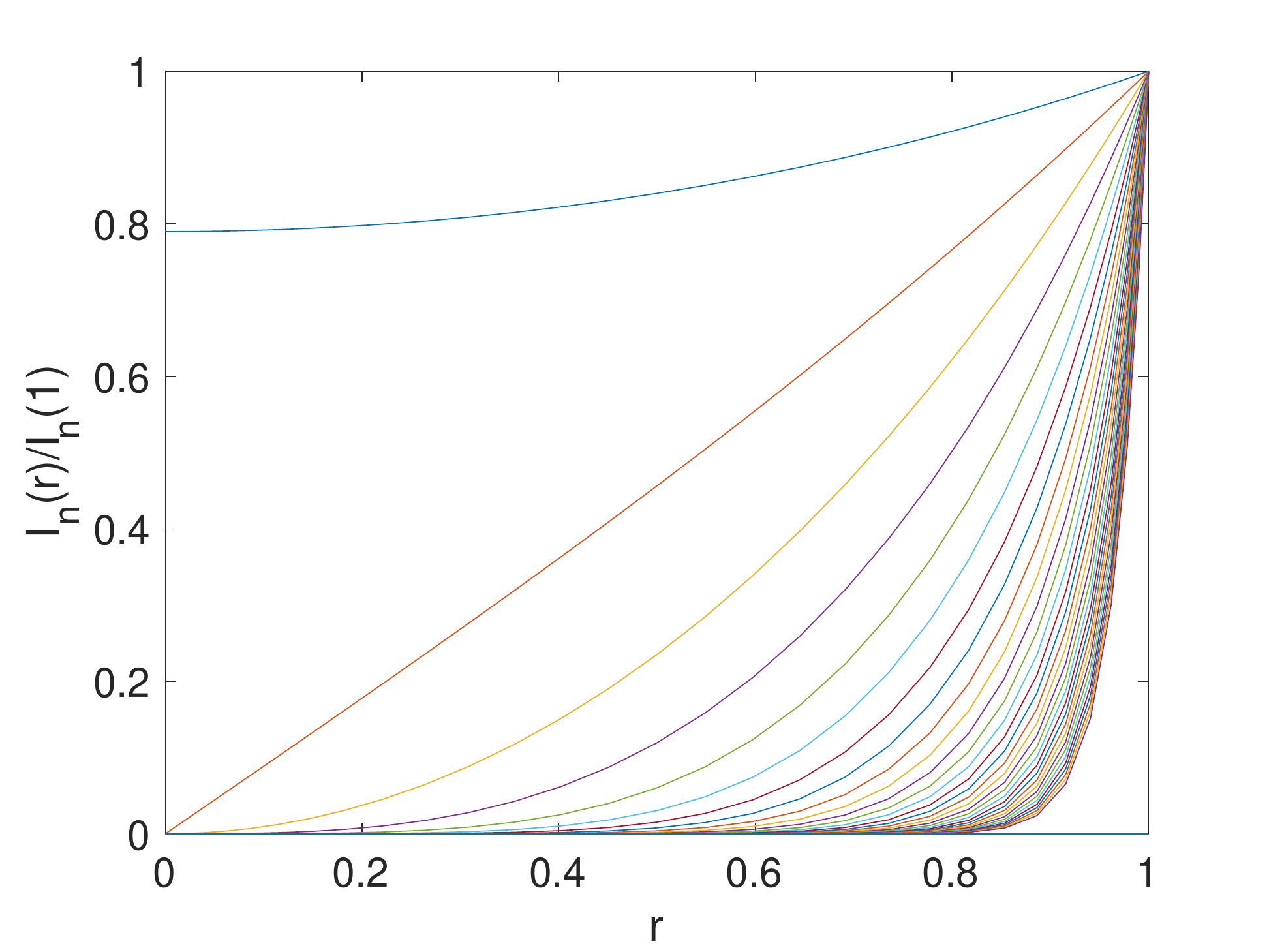}
   \includegraphics[width=0.49\textwidth]{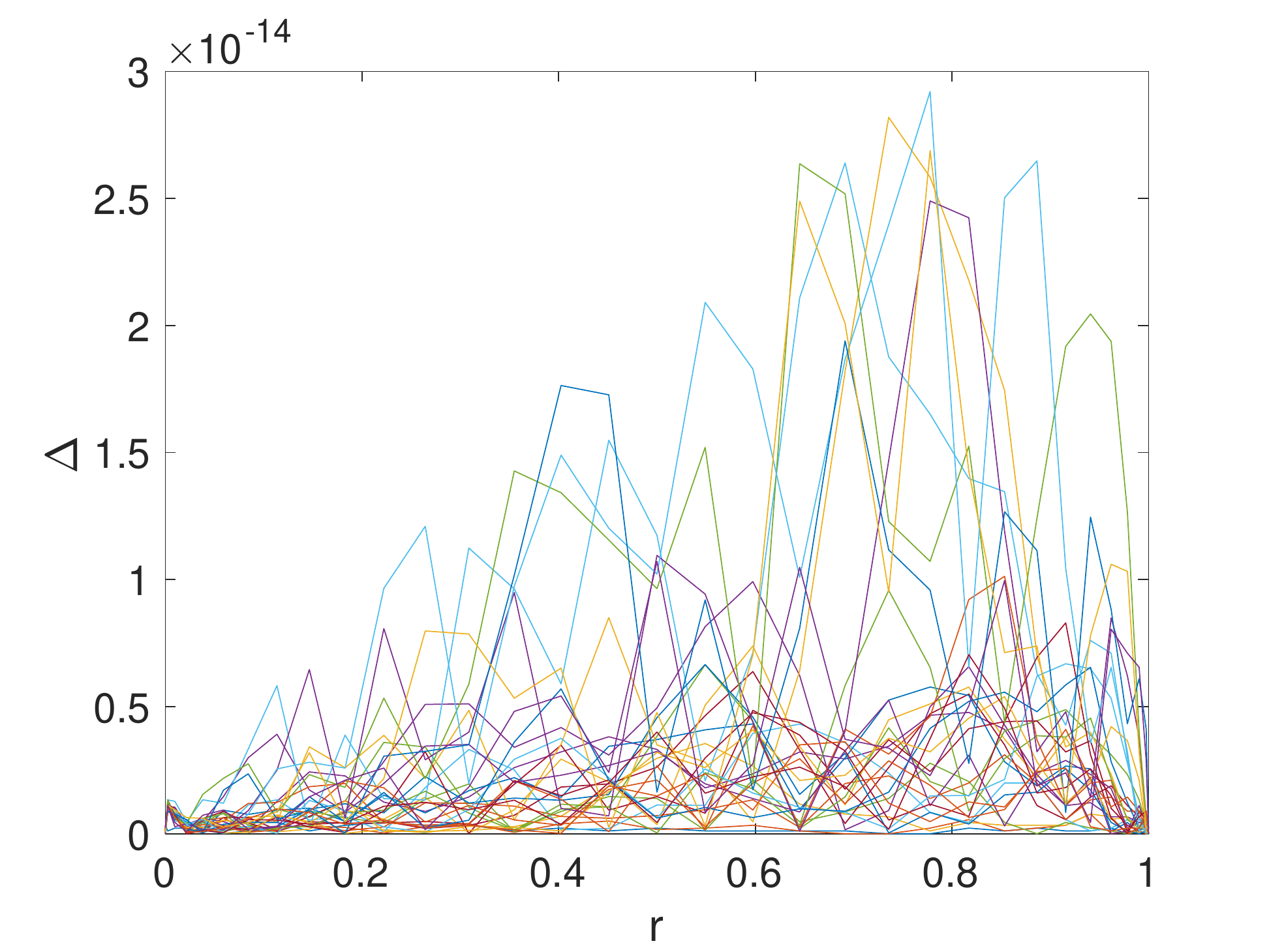}
 \caption{On the left the modified Bessel functions $I_{n}$, 
 $n=0,1,\ldots,15$ for $r\leq 1$ 
 normalized to 1 for $r=1$ (the higher $n$, the closer the 
 corresponding function is to the vertical axis on the right), and 
 the modulus of the difference between numerical and exact solution 
 for $\psi_{1}^{(j)}$, $j=1,\ldots,N_{\phi}/2$, on the right.}
 \label{besselfig}
\end{figure}

To test the above approach, we consider the example of section 
\ref{example} with $N_{r}=32$ and $N_{\phi}=64$. Relations (\ref{psi11}) and (\ref{psi22}) imply for 
$j=1,\ldots,N_{\phi}/2$ 
\begin{equation}
    \psi^{(j)}_{1}=\frac{I_{j-1}(r)}{I_{j-1}(1)}\mathrm{e}^{\mathrm{i}(j-1)\phi},\quad 
    \psi^{(j)}_{2}=\frac{I_{j}(r)}{I_{j-1}(1)}\mathrm{e}^{\mathrm{i}j\phi},
    \label{j1}
\end{equation}
and for $j=N_{\phi}/2+2,\ldots,N_{\phi}$ 
\begin{equation}
    \psi^{(j)}_{1}=\frac{I_{\tilde{j}-1}(r)}{I_{\tilde{j}}(1)}\mathrm{e}^{\mathrm{i}\tilde{j}-1\phi},\quad 
    \psi^{(j)}_{2}=\frac{I_{\tilde{j}}(r)}{I_{\tilde{j}}(1)}\mathrm{e}^{\mathrm{i}\tilde{j}\phi},
    \label{j2}
\end{equation}
where $\tilde{j}=0,-N_{\phi}/2+1,\ldots,-1$ 
for $j=N_{\phi}/2+2,\ldots,N_{\phi}$.
The Bessel functions are computed for $j\geq0$ via the series representation (see 
\cite{AS})
\begin{equation*}
    I_{j}(r)=\sum_{m=0}^{\infty}\left(\frac{r}{2}\right)^{2m+j}\frac{1}{m!(m+j)!}
\end{equation*}
The series is approximated via a finite sum consisting of the terms  
$m=0,\ldots,16$  of the series. This gives an approximation of the 
modified Bessel functions to the order of machine precision (which is 
here roughly $10^{-16}$ though in practice limited to a maximal 
accuracy of the order of $10^{-14}$ because of rounding errors). The 
functions $I_{n}$,  $n=0,\ldots,15$ can be seen in 
Fig.~\ref{besselfig} on the left. Note that $I_{-j}(r)=I_{j}(r)$ for 
$j\in\mathbb{N}$.

To study the numerical error, we define $\Delta := 
|\psi^{(j)}_{i,exact}-\psi^{(j)}_{i,numerical}|$, $i=1,2$. This error 
can be seen for $j=1,\ldots,N_{\phi}/2$ for $\psi^{(j)}_{1}$ on the 
right of Fig.~\ref{besselfig}. It is of 
the order of $10^{-14}$ on the whole disk for all values of $j$, 
i.e., of the same order of precision with which the modified Bessel 
functions are computed.

\begin{figure}[htb!]
  \includegraphics[width=0.32\textwidth]{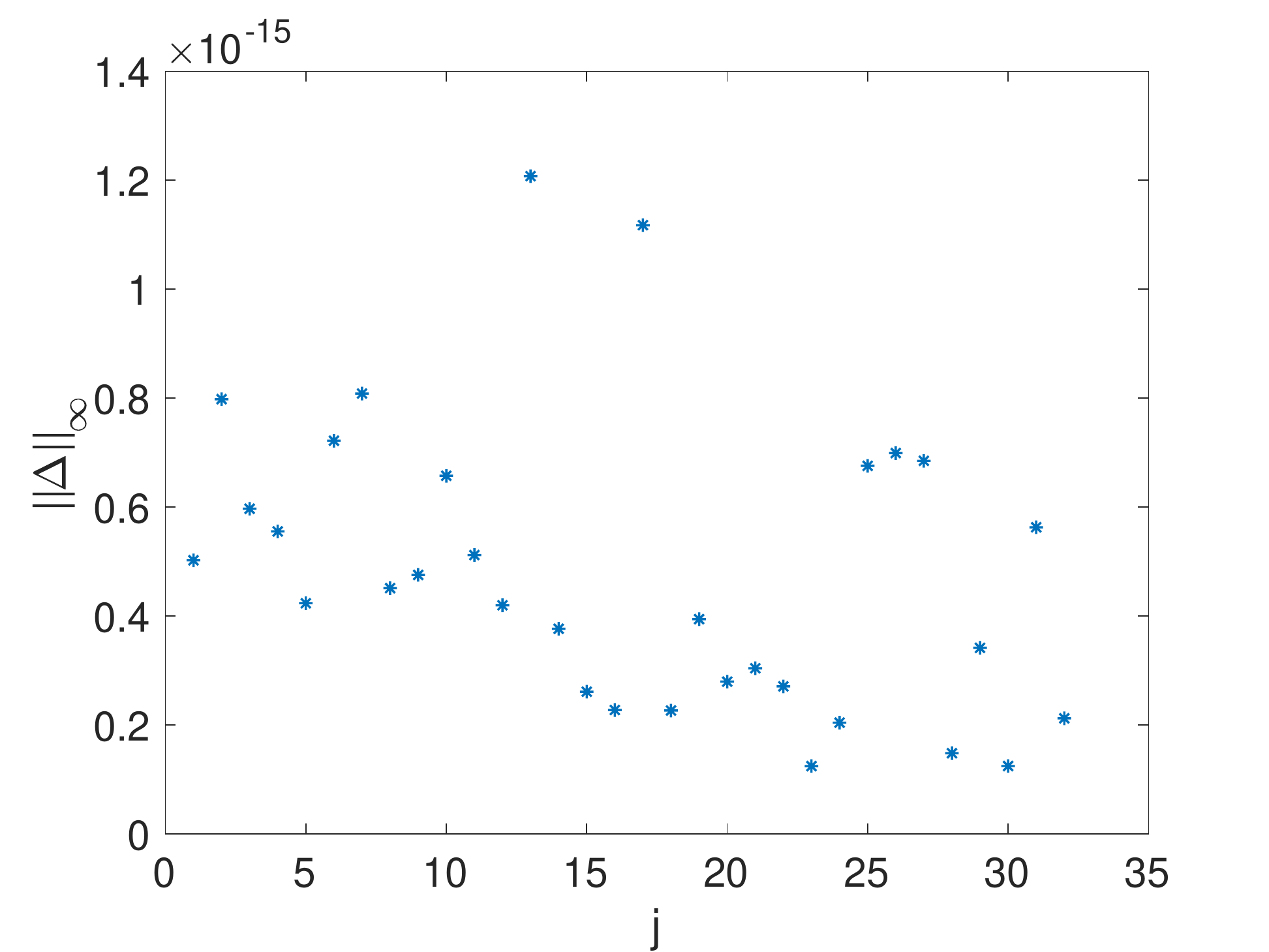}
   \includegraphics[width=0.32\textwidth]{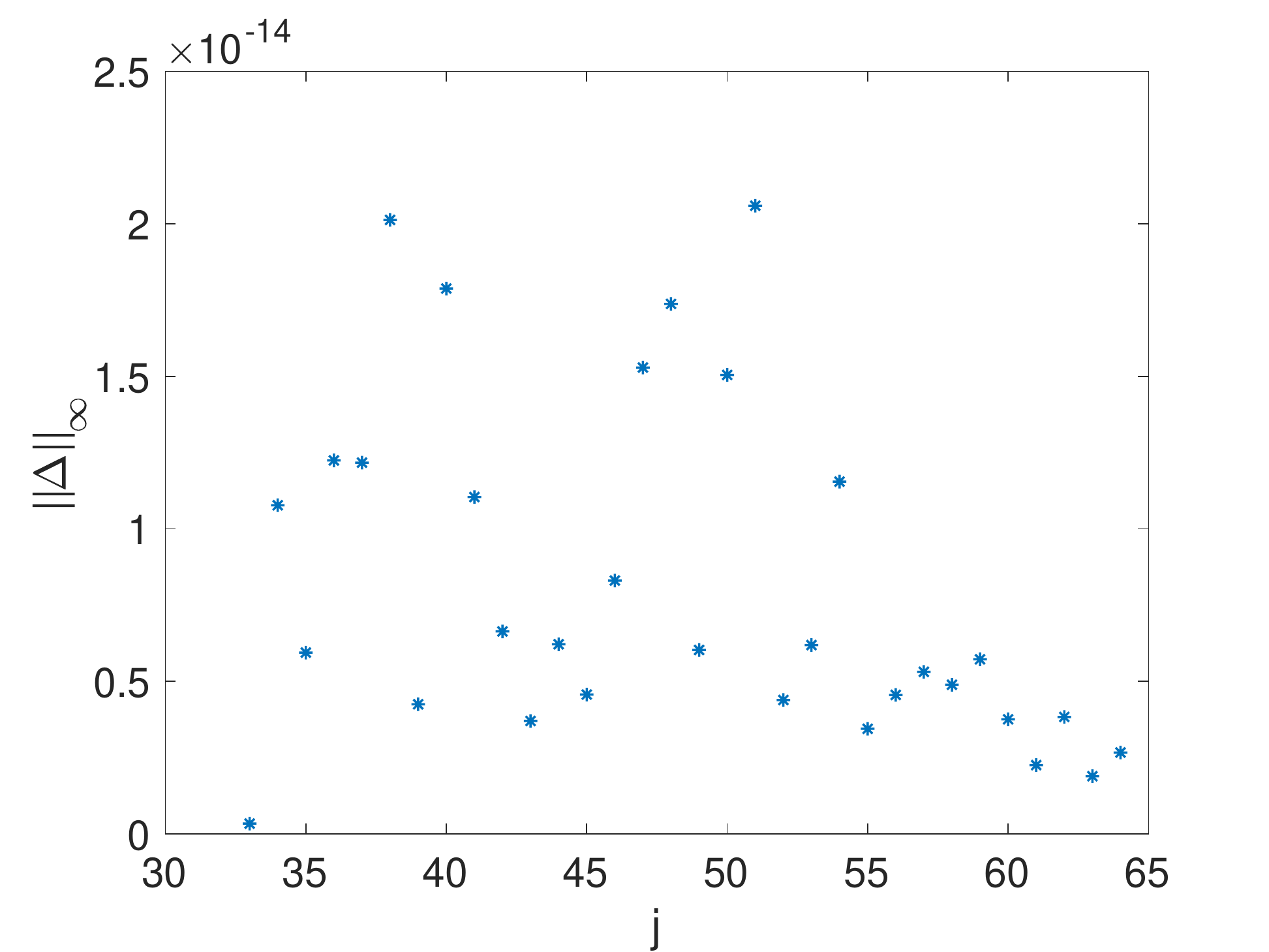}
 \includegraphics[width=0.32\textwidth]{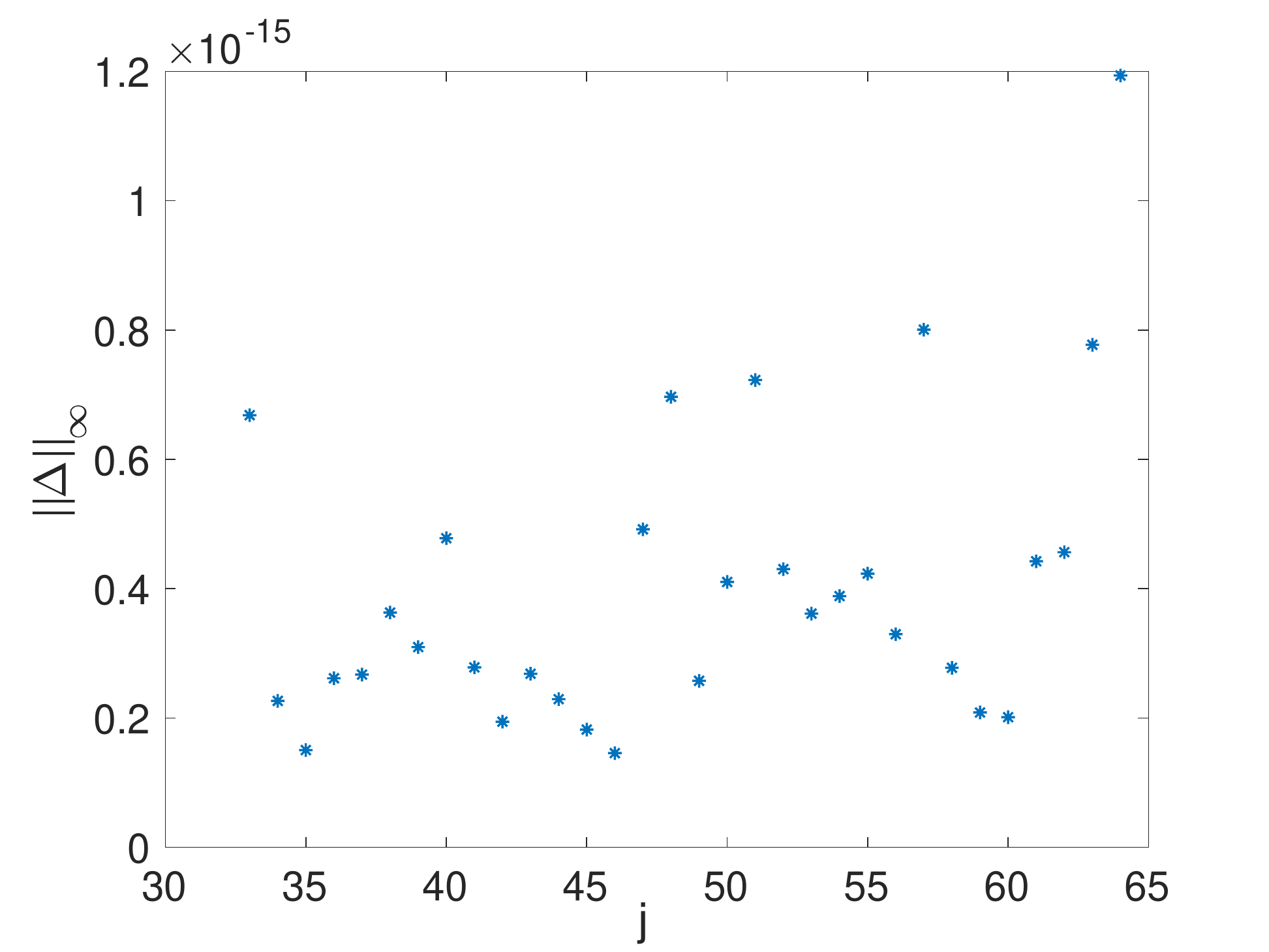}
 \caption{
$L^{\infty}$ norm of the difference between numerical and exact solution 
 for $\psi_{2}^{(j)}$, $j=1,\ldots,N_{\phi}/2$, on the left, and for 
 $j=N_{\phi}/2+1,\ldots,N_{\phi}$ and $\psi_{1}^{(j)}$ respectively 
 $\psi_{2}^{(j)}$ on the middle respectively on the right.}
 \label{besselfig2}
\end{figure}

The same level of accuracy is obtained for the functions 
$\psi^{(j)}_{2}$ for $j=1,\ldots,N_{\phi}/2$ as can be seen in 
Fig.~\ref{besselfig2} on the left. The corresponding errors for 
$j=N_{\phi}/2+1,\ldots,N_{\phi}$ are given for $\psi_{1}^{(j)}$ in 
the middle and  $\psi_{2}^{(j)}$ on the right of 
Fig.~\ref{besselfig2}. The errors are in all cases of the order of 
$10^{-14}$ or better. 

\begin{figure}[htb!]
  \includegraphics[width=0.49\textwidth]{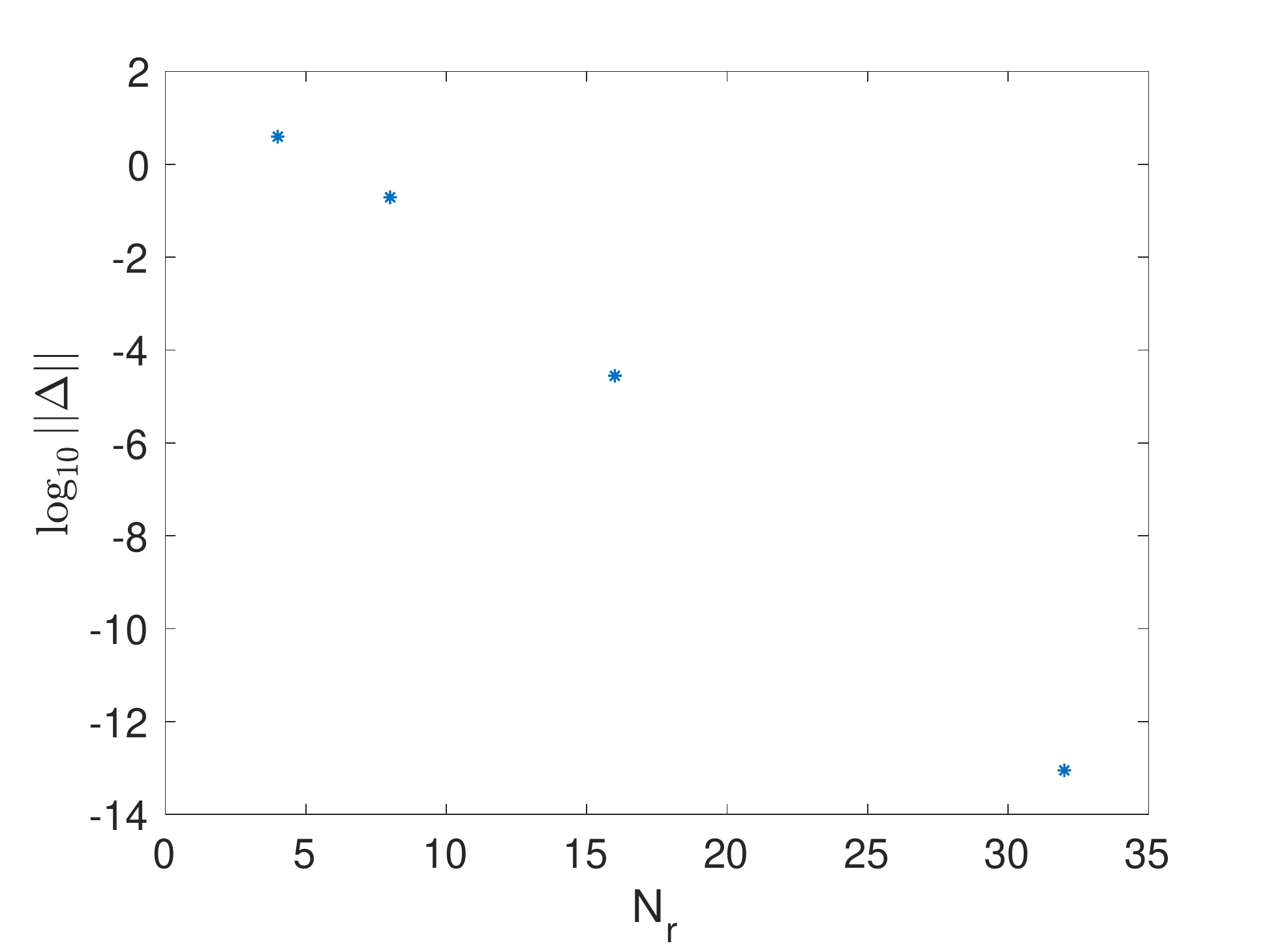}
   \includegraphics[width=0.49\textwidth]{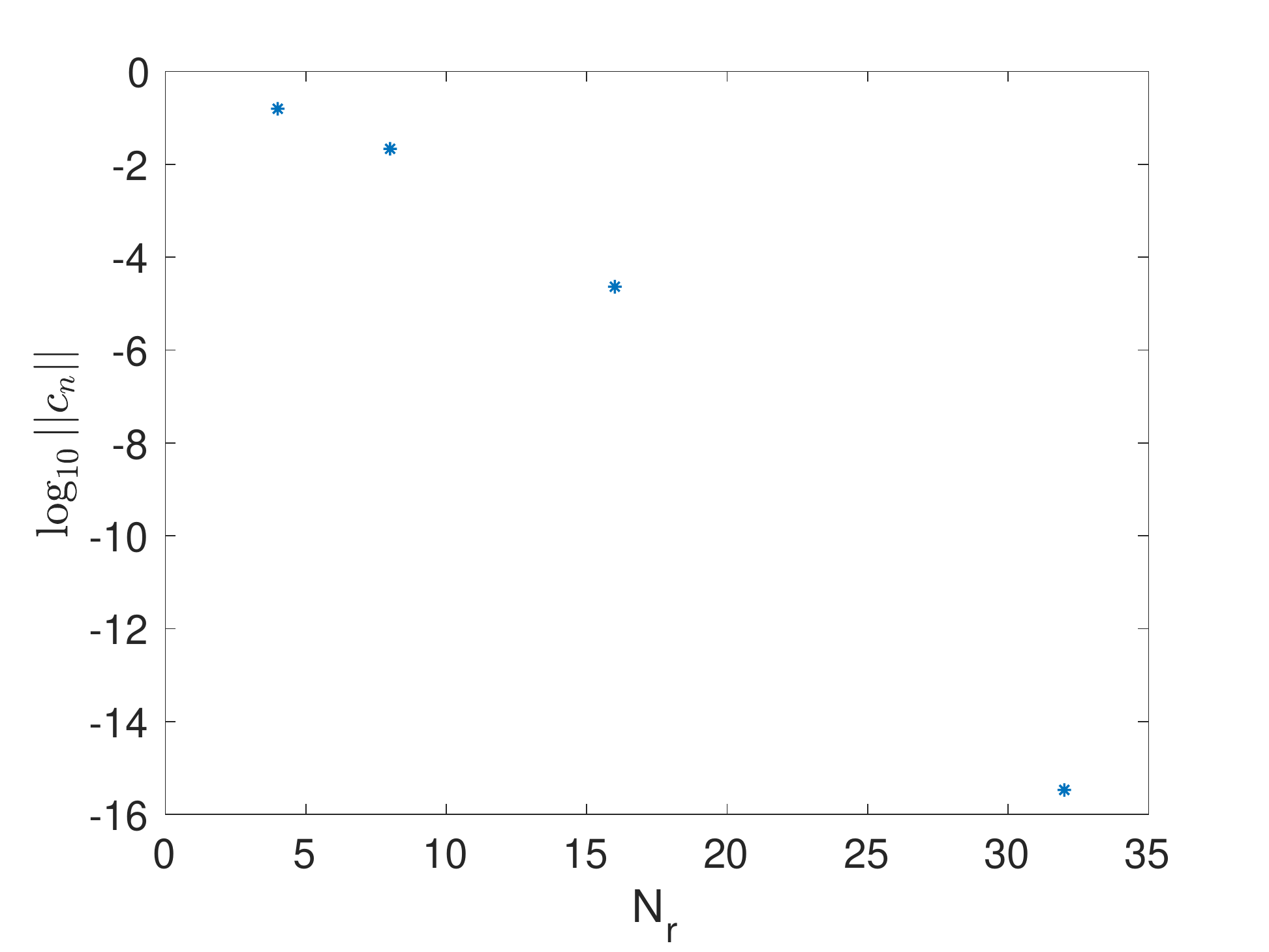}
 \caption{
 Maximum of the errors shown in Fig.~\ref{besselfig2} in dependence 
 of the number of Chebyshev modes on the left, and the maximum of the 
 highest Chebyshev coefficients for the solutions of the right.}
 \label{errorNr}
\end{figure}

An often used, but rarely in a specific context proven result in the theory of spectral methods 
is that the highest spectral coefficients indicate the numerical 
error, see for instance \cite{CC60} where it was shown that the 
error in the Clenshaw-Curtis integration is controlled by the highest 
3 spectral coefficients (note that this can be different in the case 
of equations ranging over several orders of magnitude or for 
highly ill-conditioned problems, see for instance \cite{CFKSV}). We test this behavior in the present 
context as shown in Fig.~\ref{errorNr}: on the left of the figure one 
can see that the maximum of the errors in Fig.~\ref{besselfig2} 
decreases exponentially with the number $N_{r}$ of the Chebyshev 
polynomials (the number of Fourier modes is always is always 64; as 
discussed above, due to the normalization of the solution for each 
Fourier mode, no decrease in the Fourier index can be expected) and 
the maximum of the highest Chebyshev coefficient on the right. It can 
be seen that the spectral coefficients are as expected a valid 
indicator of the numerical error and will be used in this way in the 
paper. 

\begin{figure}[htb!]
  \includegraphics[width=0.49\textwidth]{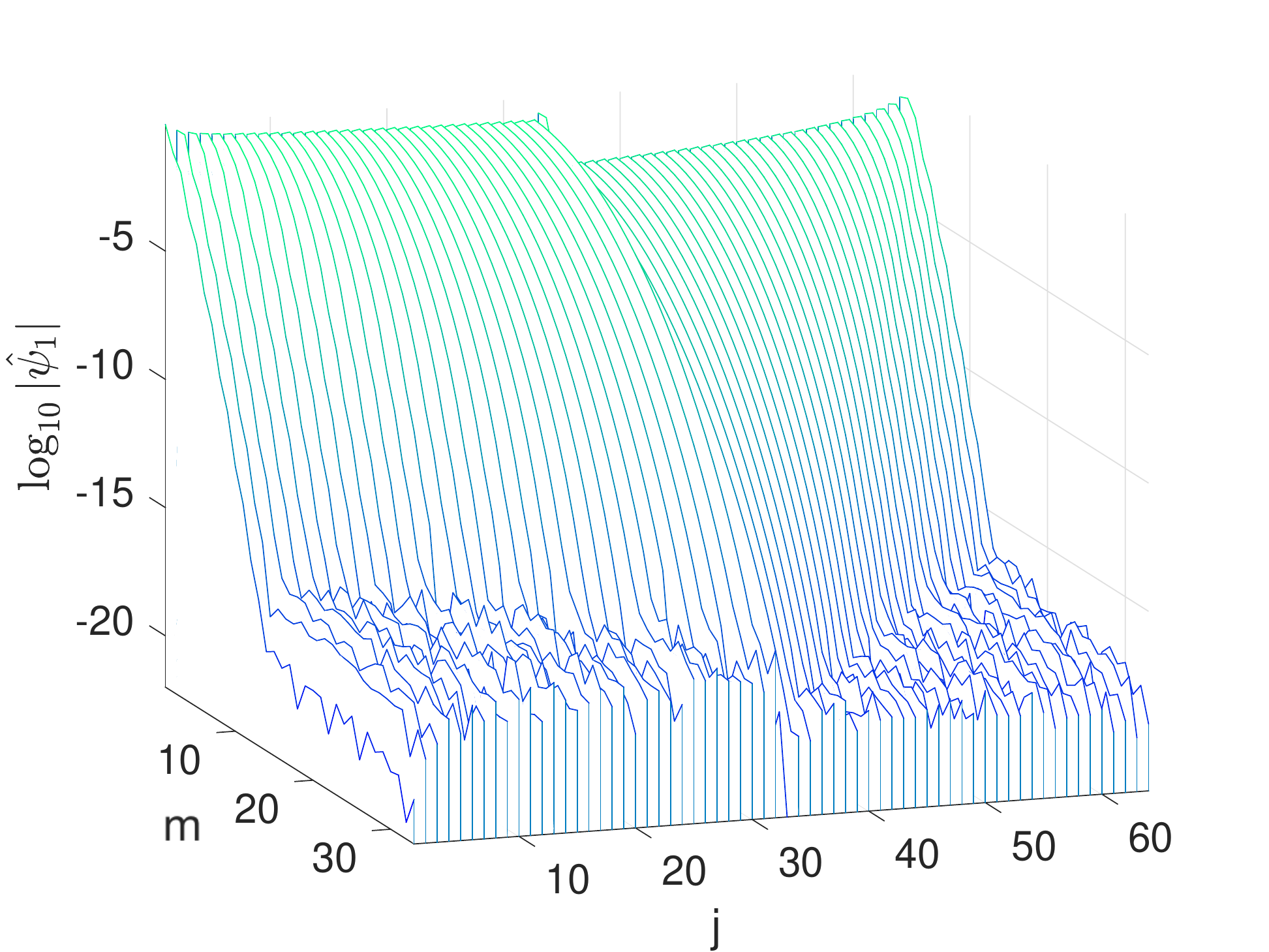}
   \includegraphics[width=0.49\textwidth]{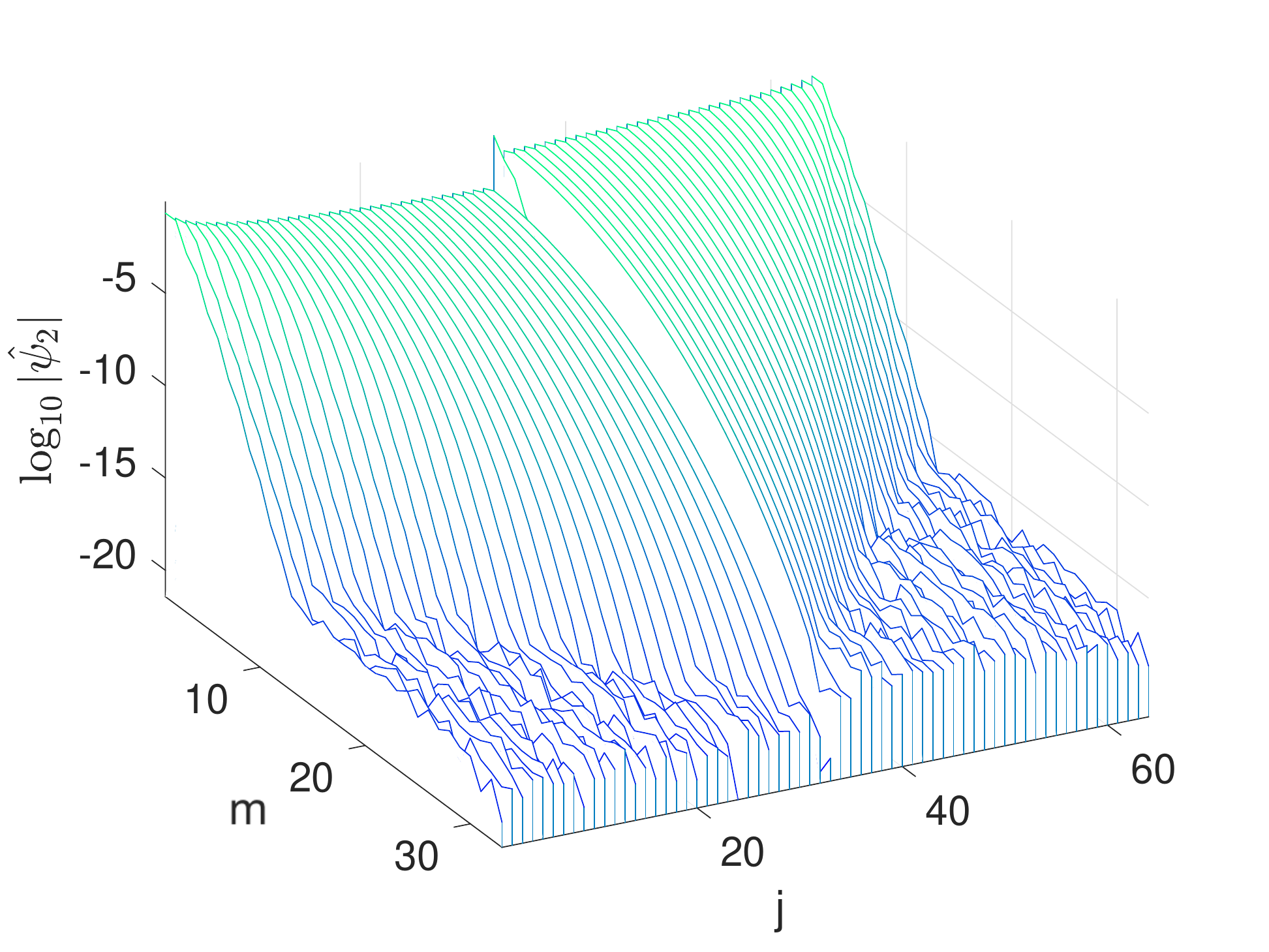}
 \caption{
 Modulus of the Chebyshev coefficients 
 $\hat{\psi}_{1}^{(j)}$ (left) and  $\hat{\psi}_{1}^{(j)}$ (right) of the functions 
 given in (\ref{j1}) and (\ref{j2}).}
 \label{besselcheb}
\end{figure}

We present the Chebyshev coefficients of the functions (\ref{j1}) and 
(\ref{j2}) respectively on the left respectively on the right of 
Fig.~\ref{besselcheb}. Evidently the coefficients decrease 
to the order of the rounding error as expected. The apparent 
discontinuity at $N_{\phi}/2$ is due to the different normalization 
conditions for $j\leq N_{\phi}/2$ (\ref{abcond1}) and for 
$j>N_{\phi}/2$ (\ref{abcond2}). It is clear that increasing the 
number of Fourier modes $N_{\phi}$ would make it necessary to choose 
a larger $N_{r}$ if spatial resolution is to be maintained. 
Since there is only 
one non-trivial Fourier coefficient per function in (\ref{j1}) and 
(\ref{j2}), we do not consider the Fourier dependence here. In the 
general case, this would be, however, necessary. 

The decrease of the 
coefficients with increasing Chebyshev and Fourier index is a test 
for the numerical approach in case there is no analytical test available. Note that the functions $\psi^{(j)}_{1,2}$ form a 
system of fundamental solutions and any solution can be expanded in terms of these functions. 
There is no decrease of the functions $\psi^{(j)}_{1,2}$ 
with the index $j$ though it is related to a Fourier index because of 
the conditions  (\ref{abcond1}) and  (\ref{abcond2}). However, we expect 
that the CGO solutions (\ref{general}) built from the functions
$\psi^{(j)}_{1,2}$ and satisfying the conditions (\ref{cnfourier}) 
and (\ref{dnfourier}) have spectral coefficients decreasing 
exponentially with Chebyshev and Fourier index as can be seen in the 
next section. This implies also that the $\gamma_{n}$, 
$n=1,\ldots,N_{\phi}$ in (\ref{general}) must decrease 
exponentially. 

\section{Complex geometric optics solutions}
In this section we use the fundamental solution constructed in the 
previous section to identify CGO solutions satisfying 
(\ref{abcond1}) and (\ref{abcond2}). To this end we determine the 
constants $\gamma_{n}$, $n=1,\ldots,N_{\phi}$ in (\ref{general}) for a
given $k$ via a linear system of equations. We consider the example 
of subsection \ref{example} to illustrate the approach. 
%

To satisfy the conditions (\ref{abcond1}) and (\ref{abcond2}), 
we consider the functions $\psi^{(j)}_{1,2}$ for $r=1$ in physical space, 
i.e., after an FFT in the variable $\phi$ for fixed $r$ for all 
$j=1,\ldots,N_{\phi}$. The coefficients $c_{n}$ (\ref{cnfourier}) and 
$d_{n}$ (\ref{dnfourier}) are computed for each $j$ from the products  $\psi^{(j)}_{1}\exp(-k\mathrm{e}^{\mathrm{i\phi}})$  and $\psi^{(j)}_{2}\exp(k\mathrm{e}^{\mathrm{-i\phi}})$ respectively 
via an FFT as per (\ref{cnfourier}) and (\ref{dnfourier}). Because of the fast algorithm this is efficient 
though only the non-negative indices $n$ are needed in the 
first case in the conditions $c_{0}=1$ and $c_{n}=0$, $n>0$, and the 
non-positive indices in the conditions $d_{n}=0$, $n\leq 0$. Notice, however, 
that the coefficient $d_{1}$ gives the reflection coefficient via 
(\ref{Rd}) which is 
therefore computed at the same time as the conditions on the 
$\gamma_{j}$. 

\begin{figure}[htb!]
  \includegraphics[width=0.49\textwidth]{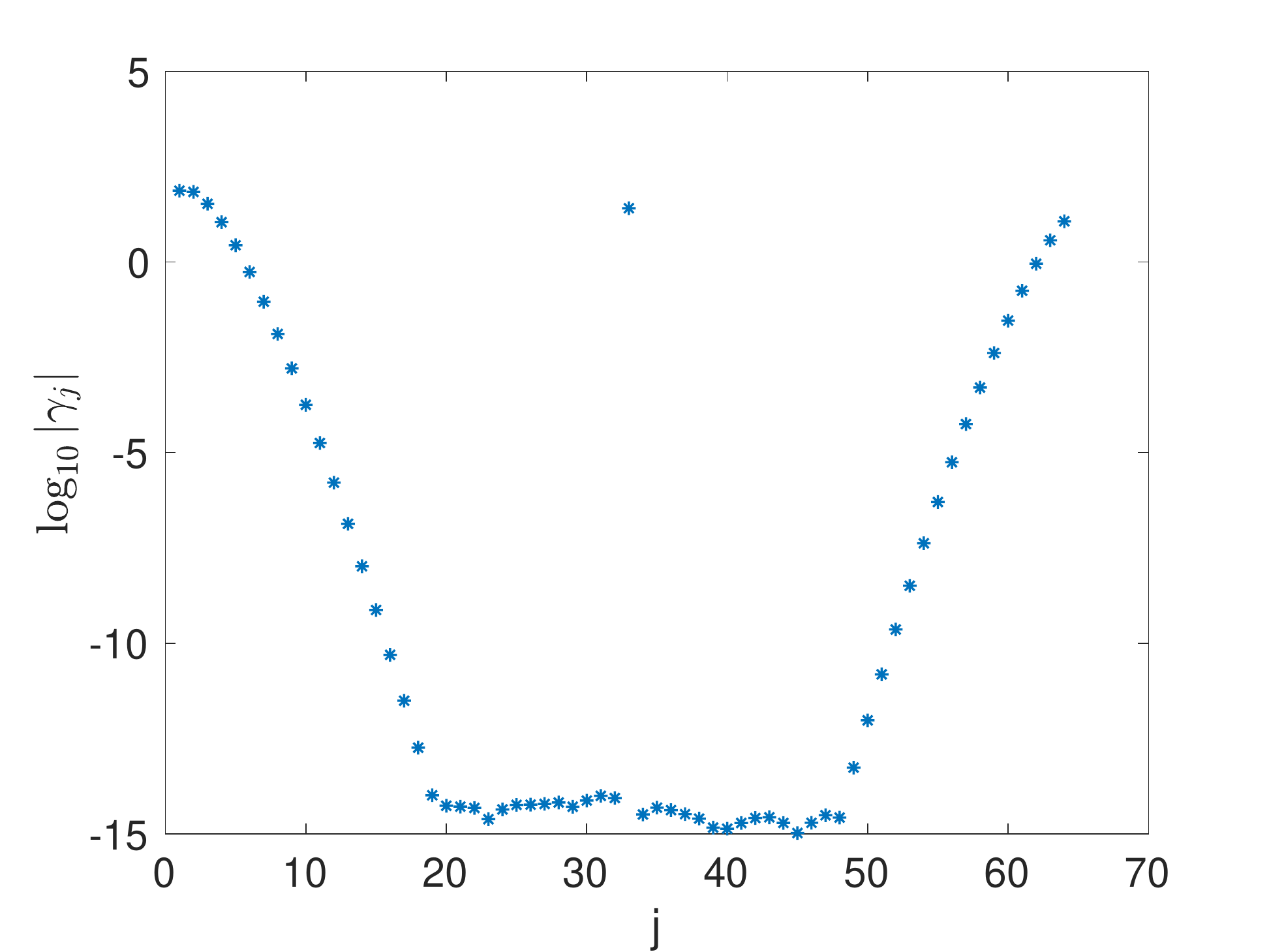}
  \includegraphics[width=0.49\textwidth]{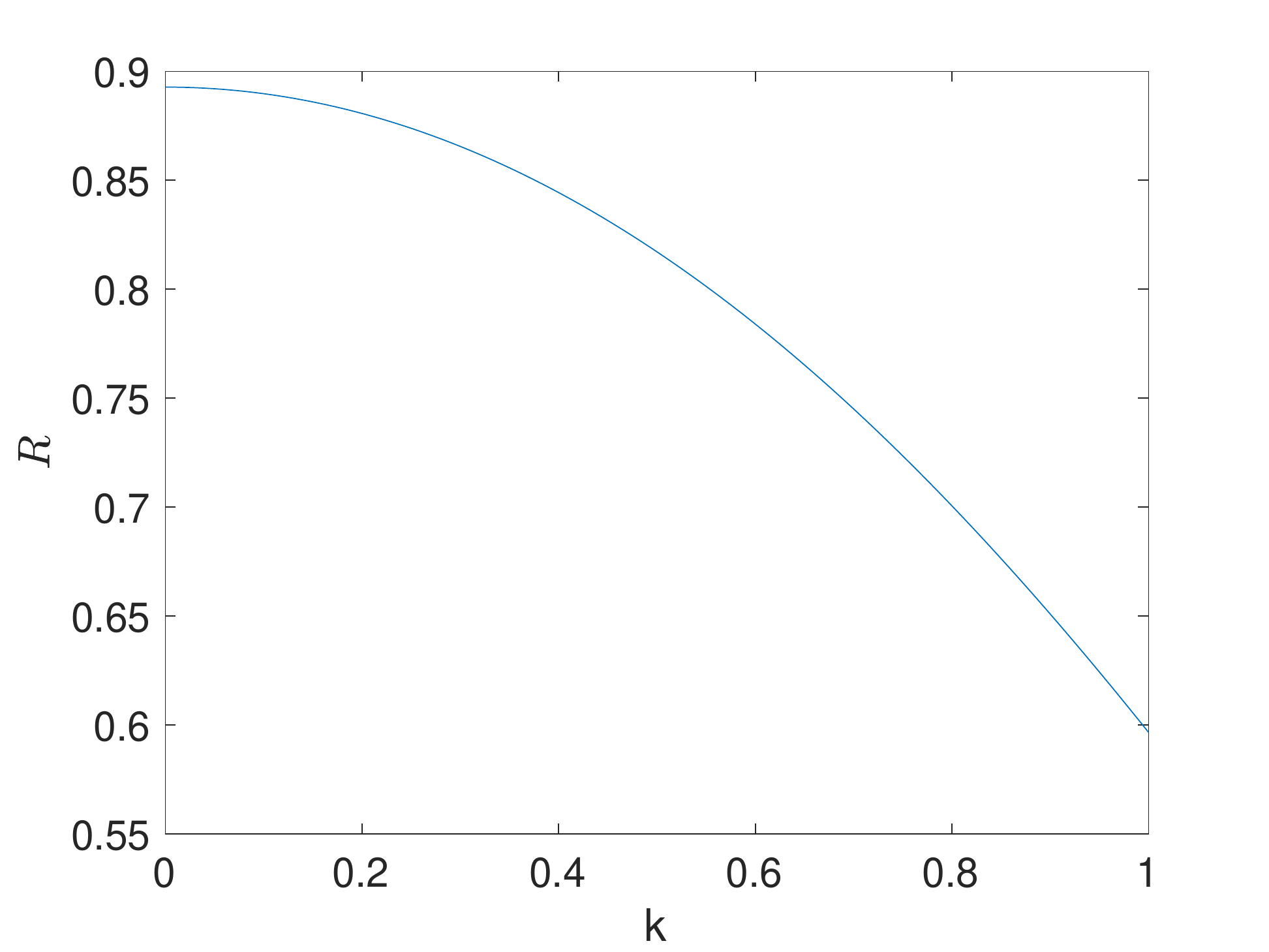}
 \caption{The coefficients $\gamma_{j}$, $j=1,\ldots,N_{\phi}$ of 
 (\ref{general}) for 
 the example $q=1$ on the disk and $k=1$ on the left, and the 
 reflection coefficient in dependence of $k$ on the right.}
 \label{gammak1}
\end{figure}

The conditions on $c_{n}$ and $d_{n}$ lead to an $N_{\phi}$ 
dimensional linear system of equations for the $\gamma_{j}$, 
$j=1,\ldots,N_{\phi}$ in (\ref{general}).  For the example of 
Subsection \ref{example}, the coefficients $\gamma_{j}$, $j=1,\ldots,N_{\phi}$
can be seen in Fig.~\ref{gammak1} on the left.  As expected, the coefficients 
 for $j\sim N_{\phi}/2$ are of the order of the rounding error (note 
though that $\gamma_{N_{\phi}/2+1}$ corresponds to the case $b_{0}=1$ and is 
thus as expected of order 1). 

\begin{figure}[htb!]
  \includegraphics[width=0.49\textwidth]{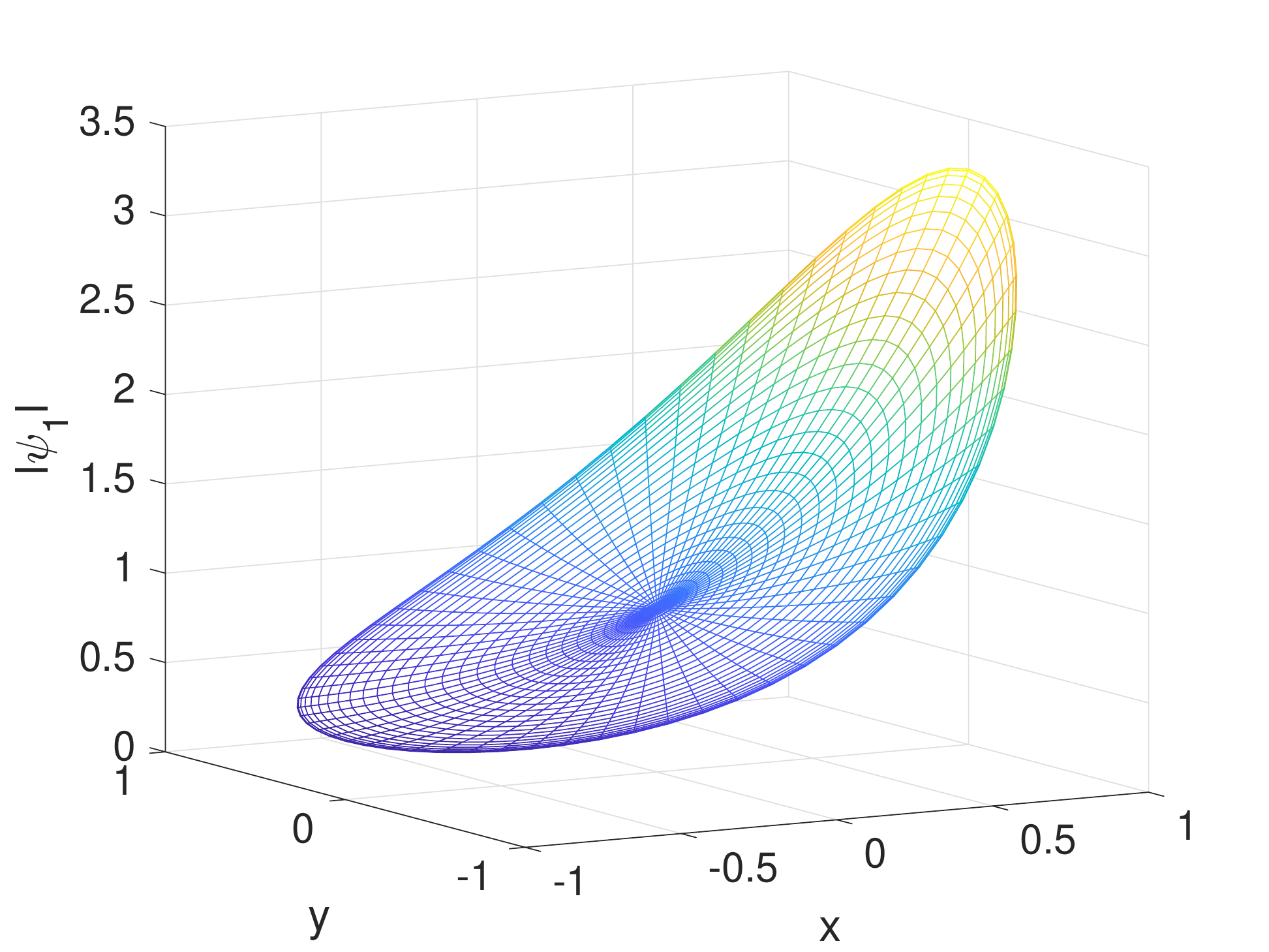}
   \includegraphics[width=0.49\textwidth]{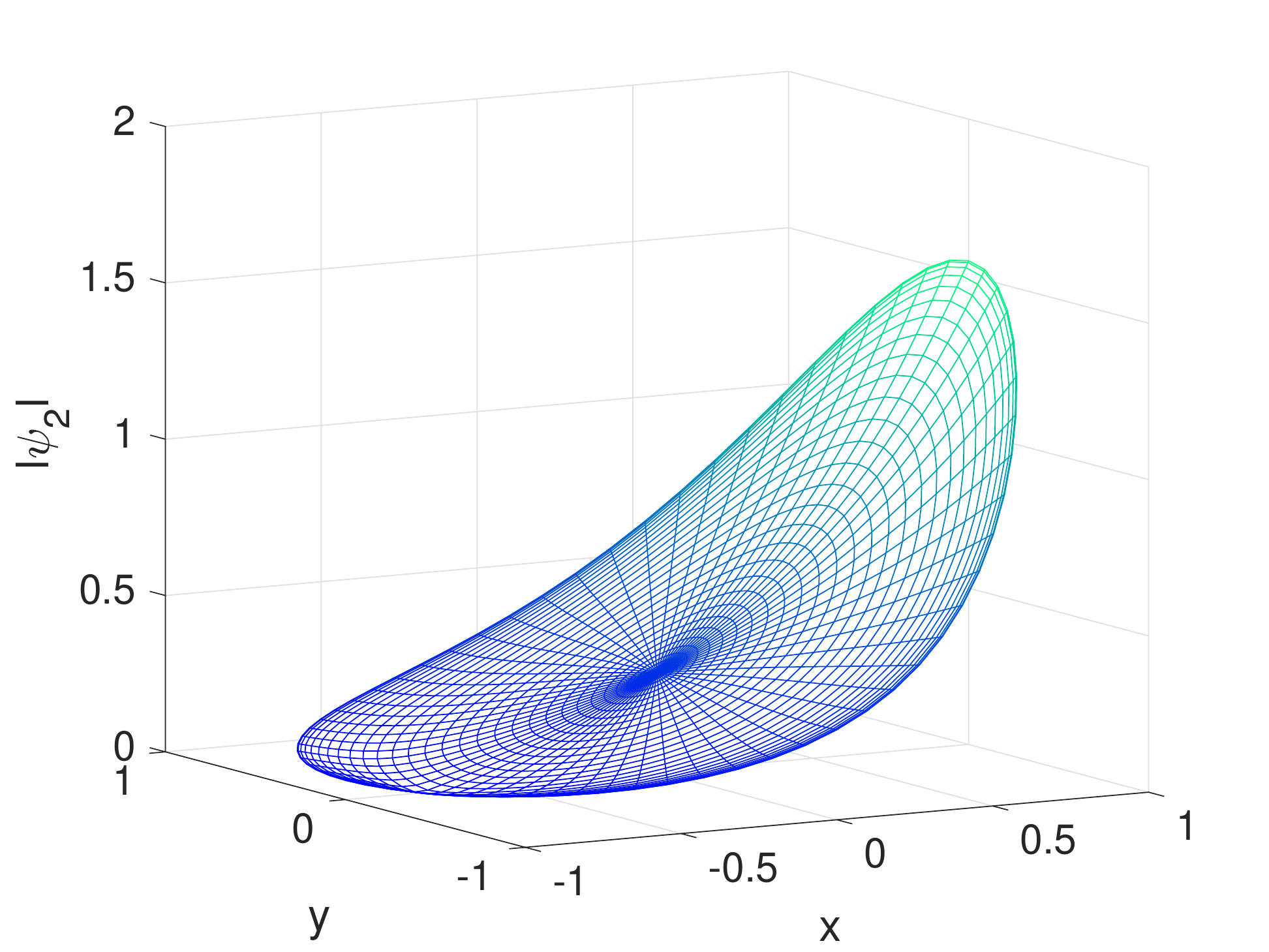}
 \caption{CGO solutions to (\ref{dbarpsi}) and (\ref{dbarpsiasym}) 
 for $q=1$ on the unit disk and vanishing elsewhere for $k=1$, on the 
 left $\psi_{1}$, on the right $\psi_{2}$.}
 \label{Phik1}
\end{figure}

The corresponding solutions (\ref{general}) can be seen in 
Fig.~\ref{Phik1} on the disk. Note that the solutions are diverging 
for $|z|\to\infty$. 

\begin{figure}[htb!]
  \includegraphics[width=0.49\textwidth]{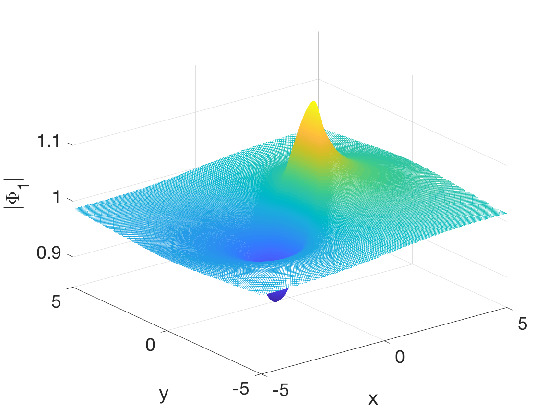}
   \includegraphics[width=0.49\textwidth]{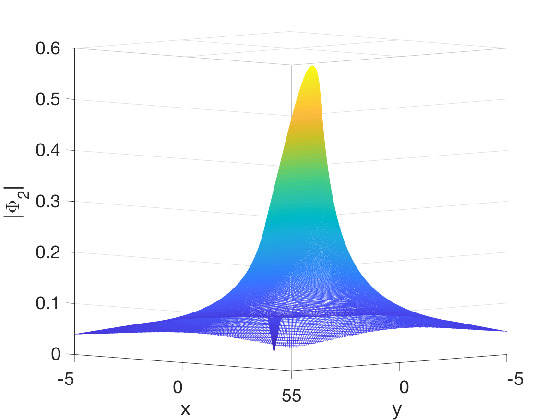}
 \caption{CGO solutions (\ref{Phi}) to (\ref{dbarpsi}) and (\ref{dbarpsiasym}) 
 for $q=1$ on the unit disk and vanishing elsewhere for $k=1$, on the 
 left $\Phi_{1}$, on the right $\Phi_{2}$.}
 \label{Phik1exp}
\end{figure}

On the other hand, the functions $\Phi_{1}$ and 
$\Phi_{2}$ (\ref{Phi}) as shown in Fig.~\ref{Phik1exp} 
are bounded for all $z\in \mathbb{C}$. 

\begin{figure}[htb!]
  \includegraphics[width=0.49\textwidth]{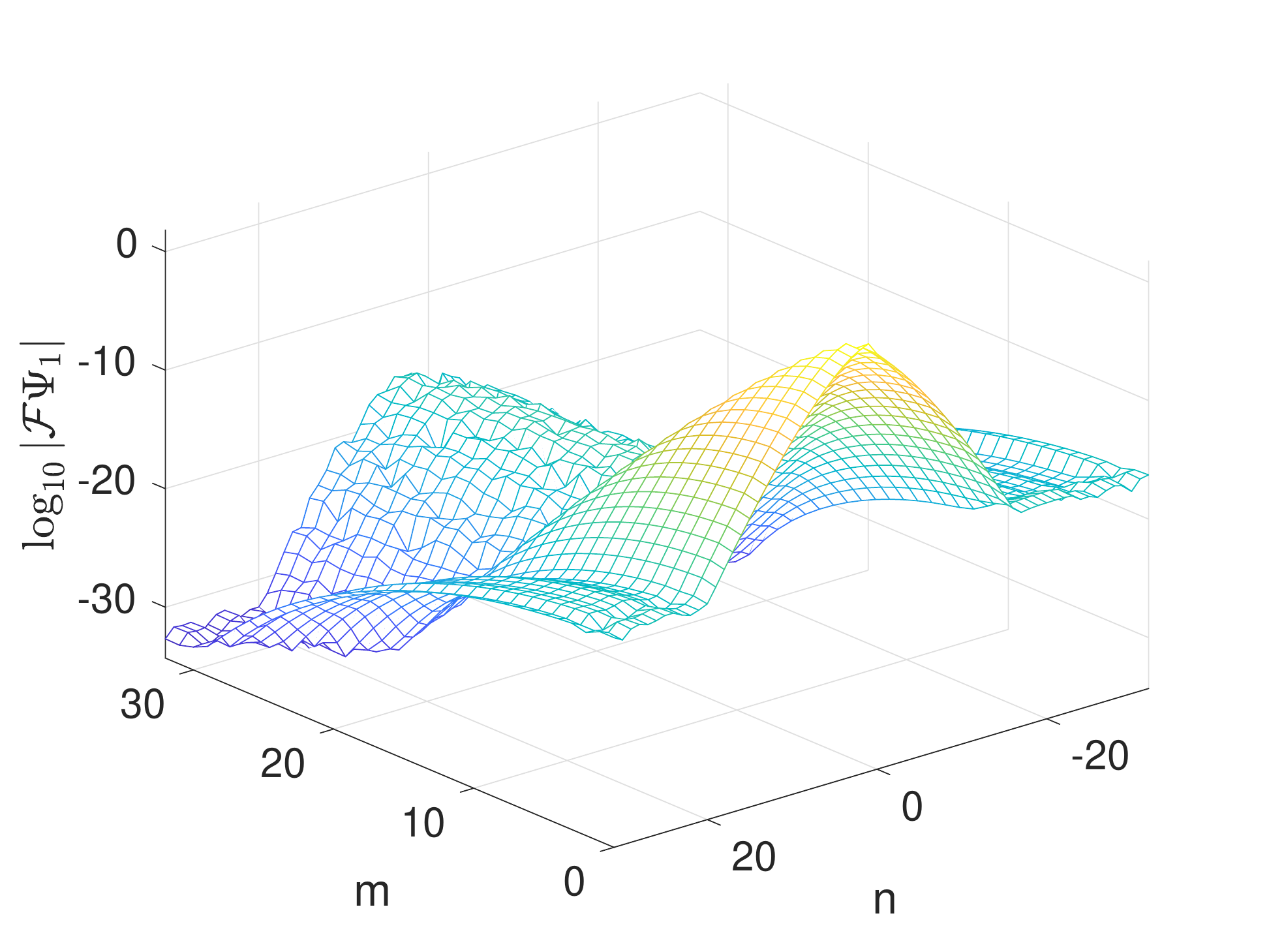}
   \includegraphics[width=0.49\textwidth]{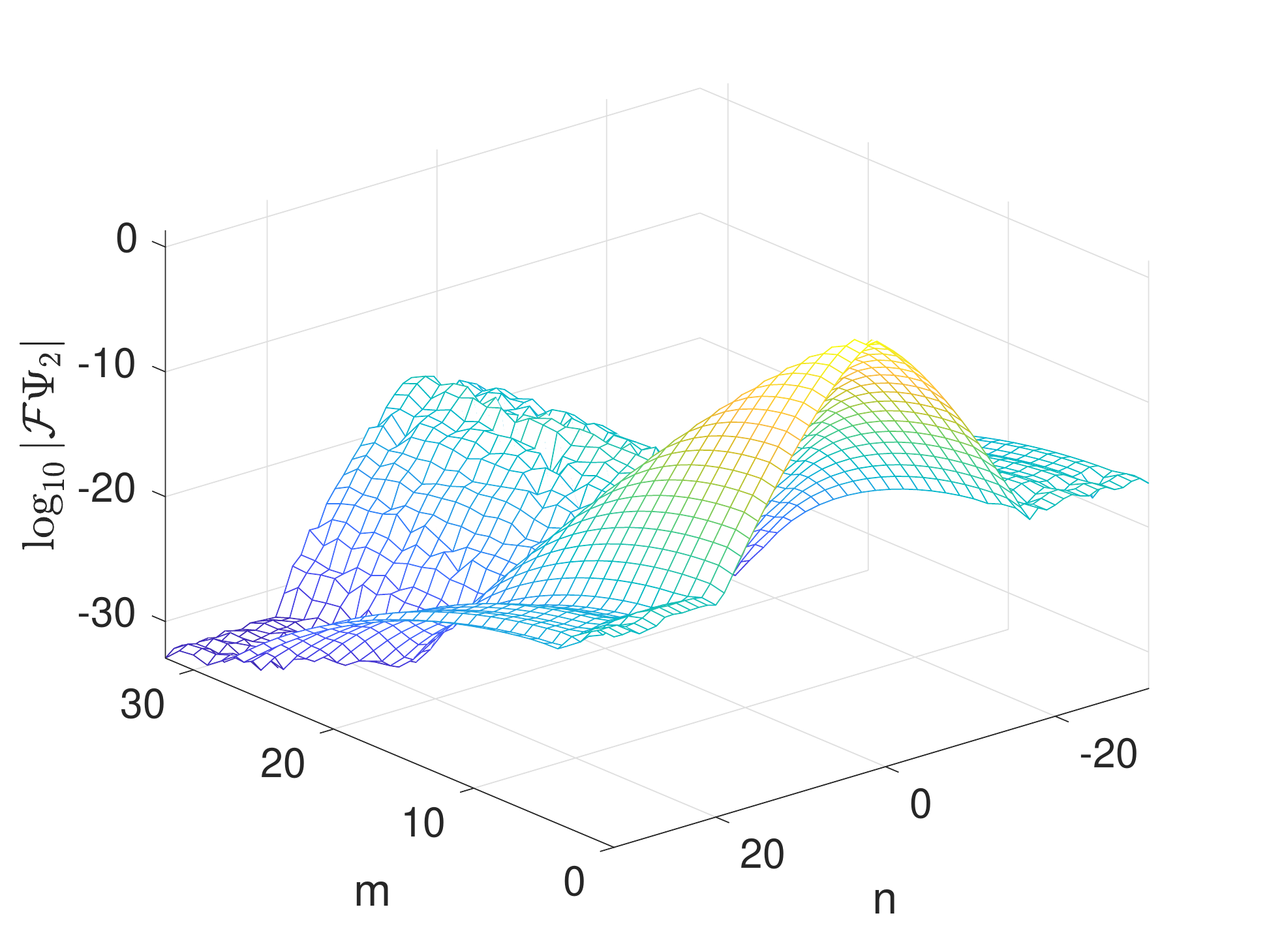}
 \caption{
 Spectral coefficients of the CGO solutions to (\ref{dbarpsi}) and (\ref{dbarpsiasym}) 
 for $q=1$ on the unit disk and vanishing elsewhere for $k=1$, on the 
 left for $\psi_{1}$, on the right for $\psi_{2}$.}
 \label{Phik1coeff}
\end{figure}

The spectral coefficients for the solutions to (\ref{dbarpsi}) and (\ref{dbarpsiasym})  are 
shown in Fig.~\ref{Phik1coeff}. It can be seen that the coefficients 
decrease to machine precision both in the Chebyshev and in 
the Fourier index as already expected from the coefficients $\gamma_{j}$ in 
 Fig.~\ref{gammak1}. This implies that the solution is resolved to the 
order of machine precision both in $r$ and in $\phi$. 

As already mentioned, the reflection coefficient (\ref{Rd}) is given 
by the scalar product of the $\gamma_{j}$ and the coefficients 
$d_{1}^{(j)}$ computed according to (\ref{dnfourier}) for each 
$\psi_{2}^{(j)}$, $j=1,\ldots,N_{\phi}$ in the determination of the 
former. Thus it is computed at negligible cost. We show the 
reflection coefficient for   $k\leq 1$ for the example $q=1$ on the 
right of Fig.~\ref{gammak1}. Note that the reflection coefficient 
depends only on $|k|$ because of the radial symmetry of $q(z)$, and 
that it is real in the studied example. 

\begin{remark}
It is of course possible to numerically implement the CGO conditions 
(\ref{cnfourier}) and  (\ref{dnfourier}) for given $k$ instead of 
the conditions (\ref{abcond1}) and (\ref{abcond2}) in the same way as 
the latter. If one is only interested in the solution for one value 
of $k$, this is in fact more economic. In applications these 
solutions are  in principle needed for all values of $k\in \mathbb{C}$ 
which means for a grid of values of $k$. The computational cost in 
producing the fundamental solution is $N_{\phi}$ times larger than 
producing the CGO solution for given $k$. Satisfying the CGO 
conditions for a given fundamental solution is then a lower 
dimensional problem since it is done for $r=1$ only. Thus as soon as 
one computes the CGO solutions for considerably more than $N_{\phi}$ 
values of $k$, the method presented in this section is clearly more 
economical: the main computational cost is in the determination of the 
fundamental solution of the last section, the CGO conditions being of 
lower dimension are then essentially for free.
   
\end{remark}

\section{Iterative solution of the d-bar system}
In this section we present an iterative solution to the d-bar system 
(\ref{dbarpolphi}). The same spectral approach as in the previous 
sections is applied to approximate the derivatives. The resulting 
system is solved with a Picard iteration which was shown in 
\cite{joh} to 
converge for  more regular potentials for large $|k|$ like a geometric series in $1/|k|$.  Here we show at 
concrete examples that the iteration converges very rapidly, and this 
for all values of $|k|$ we can access. This allows to test the codes of the previous 
section and the iterative approach given here by comparing their 
respective results. An asymptotic formula for the reflection 
coefficent is presented. 

\subsection{Fixed point iteration}
The system (\ref{dbarpolphi}) will be solved with a fixed point 
iteration ($j=0,1,\ldots$)
\begin{equation}
\begin{array}{l}
    \mathrm{e}^{\mathrm{i}\phi}\left(\partial_{r}+\frac{\mathrm{i}}{r}\partial_{\phi}\right)\Phi_{1}^{(j+1)}  
    =q(r,\phi)\mathrm{e}^{\bar{k}\bar{z}-kz}\Phi_{2}^{(j)},
    \\~\\
     \mathrm{e}^{-\mathrm{i}\phi}\left(\partial_{r}-\frac{\mathrm{i}}{r}\partial_{\phi}\right)\Phi_{2}^{(j+1)}  
    =\bar{q}(r,\phi)\mathrm{e}^{kz-\bar{k}\bar{z}}\Phi_{1}^{(j)}
    \label{dbarpolphiit},
\end{array}
\end{equation}
with the initial iterates $\Phi^{(0)}_{1}=1$ and $\Phi^{(0)}_{2}=0$.

As in the previous sections, see (\ref{cheb}), we will apply a
Chebyshev spectral method in $r$ and a discrete Fourier approach 
in $\phi$, 
\begin{equation*}
\begin{array}{l}
    \Phi^{(j)}_{1}\approx\sum_{\alpha=0}^{N_{r}}\sum_{\beta=-N_{\phi}/2+1}^{
    N_{\phi}/2}a_{\beta\alpha}^{(j)}T_{\alpha}(l)\mathrm{e}^{\mathrm{i}\beta\phi}
,\\~\\
 \Phi^{(j)}_{2}\approx\sum_{\alpha=0}^{N_{r}}\sum_{\beta=-N_{\phi}/2}^{
    N_{\phi}/2-1}b_{\beta\alpha}^{(j)}T_{\alpha}(l)\mathrm{e}^{\mathrm{i}\beta\phi}
\end{array}
\end{equation*}

System (\ref{dbarpolphiit}) can thus be approximated by 
($m=0,1,\ldots,N_{r}$)
\begin{align}
    \sum_{\alpha=0}^{N_{r}}\left(D-nR\right)_{m\alpha}a_{n\alpha}^{(j+1)}&=\mathbb{F}\left(q(r,\phi)\mathrm{e}^{\bar{k}\bar{z}-kz-\mathrm{i}\phi}\Phi_{2}^{(j)}\right)_{nm},\quad n=-N_{\phi}/2,\ldots,N_{\phi}/2-1,
    \nonumber\\
        \sum_{\alpha=0}^{N_{r}}\left(D+nR\right)_{m\alpha}b_{n\alpha}^{(j+1)}&=\mathbb{F}\left(\bar{q}(r,\phi)
    \mathrm{e}^{kz-\bar{k}\bar{z}+{\mathrm{i}\phi}}\Phi_{1}^{(j)}\right)_{nm},\quad n=-N_{\phi}/2+1,\ldots,N_{\phi}/2.
    \label{dbarsysnumit}
\end{align}
The matrices $D\mp n R$ are inverted by using the boundary 
conditions (\ref{Phisasym}) which imply 
\begin{equation*}
    \sum_{\alpha=0}^{N_{r}}a_{n\alpha}^{(j)}=\delta_{n0}, \quad 
    n=0,\ldots,N_{\phi}/2-1,
\end{equation*}
and 
\begin{equation}
    \sum_{\alpha=0}^{N_{r}}b_{n\alpha}^{(j)}=0, \quad 
    n=-N_{\phi}/2+1,\ldots,0.
    \label{phi2cond}
\end{equation}
Note that this can be seen as a simplified Newton iteration since just 
the diagonal part of the Jacobian 
$$\mbox{Jac}=
\begin{pmatrix}
    \bar{\partial} & -\mathrm{e}^{\bar{k}\bar{z}-kz}q/2 \\
     -\mathrm{e}^{kz-\bar{k}\bar{z}}\bar{q}/2& \partial
\end{pmatrix}
$$ 
(the d-bar system is linear and thus equal to $\mbox{Jac}
\begin{pmatrix}
    \Phi_{1} \\
    \Phi_{2}
\end{pmatrix}
=0$) is inverted with some boundary conditions.

The right-hand sides of (\ref{dbarsysnumit}) are computed in physical 
space to avoid convolutions. Due to the fast algorithms 
 FFT and the FCT, this is very efficient. In the algorithm of 
the previous sections just one matrix $\mathcal{O}$ in (\ref{mat2}) had to be 
inverted whereas here an iterative approach is used. However, the matrix 
$\mathcal{O}$ in (\ref{mat2}) had to be determined via convolutions, 
whereas we are dealing here just with products in physical space 
and fast transforms.

\begin{remark}\label{remresolution}
    For each value of $n$ in (\ref{dbarsysnumit}), an 
    $(N_{r}+1)\times(N_{r}+1)$ matrix has to be 
    inverted which means $2N_{\phi}$ inversions of this size per 
    iteration. The advantage with respect to the approach of the previous 
    sections is that the system decouples here, just matrices of size 
    $(N_{r}+1)\times (N_{r}+1)$ and $N_{\phi}\times N_{\phi}$ have to handled. 
    This is much less demanding in terms of memory than the previous 
    $2(N_{r}+1)N_{\phi}\times 2(N_{r}+1)N_{\phi}$ matrices in (\ref{mat2})  
    (though the 
    diagonal form of the differentiation in Fourier leads to sparse 
    matrices there for which efficient algorithms are implemented in 
    Matlab). Thus a much higher resolution can be reached here 
    on the same hardware than before which is especially interesting 
    for the case of large $|k|$. In other words, the price for the 
    reduced memory requirements with respect to the method of the 
    previous sections (where one big matrix had to inverted without 
    iterating) is an iteration. The values of $N_{r}$ and $N_{\phi}$ 
    have to be essentially chosen in a way that the factor 
    $\mathrm{e}^{\bar{k}\bar{z}-kz}q$ is numerically resolved for 
    given $k$ and $q$
    on the disk.
\end{remark}

The reflection coefficient $R$ (\ref{eq:r-def}) follows in this 
approach from the function $\Phi_{2}$ for $r=1$. One gets from (\ref{phi2cond})
\begin{equation}
    R = 2\sum_{\alpha=0}^{N_{r}}\bar{b}_{1\alpha}
    \label{reflection}.
\end{equation}

\subsection{Example}
\begin{figure}[htb!]
  \includegraphics[width=0.49\textwidth]{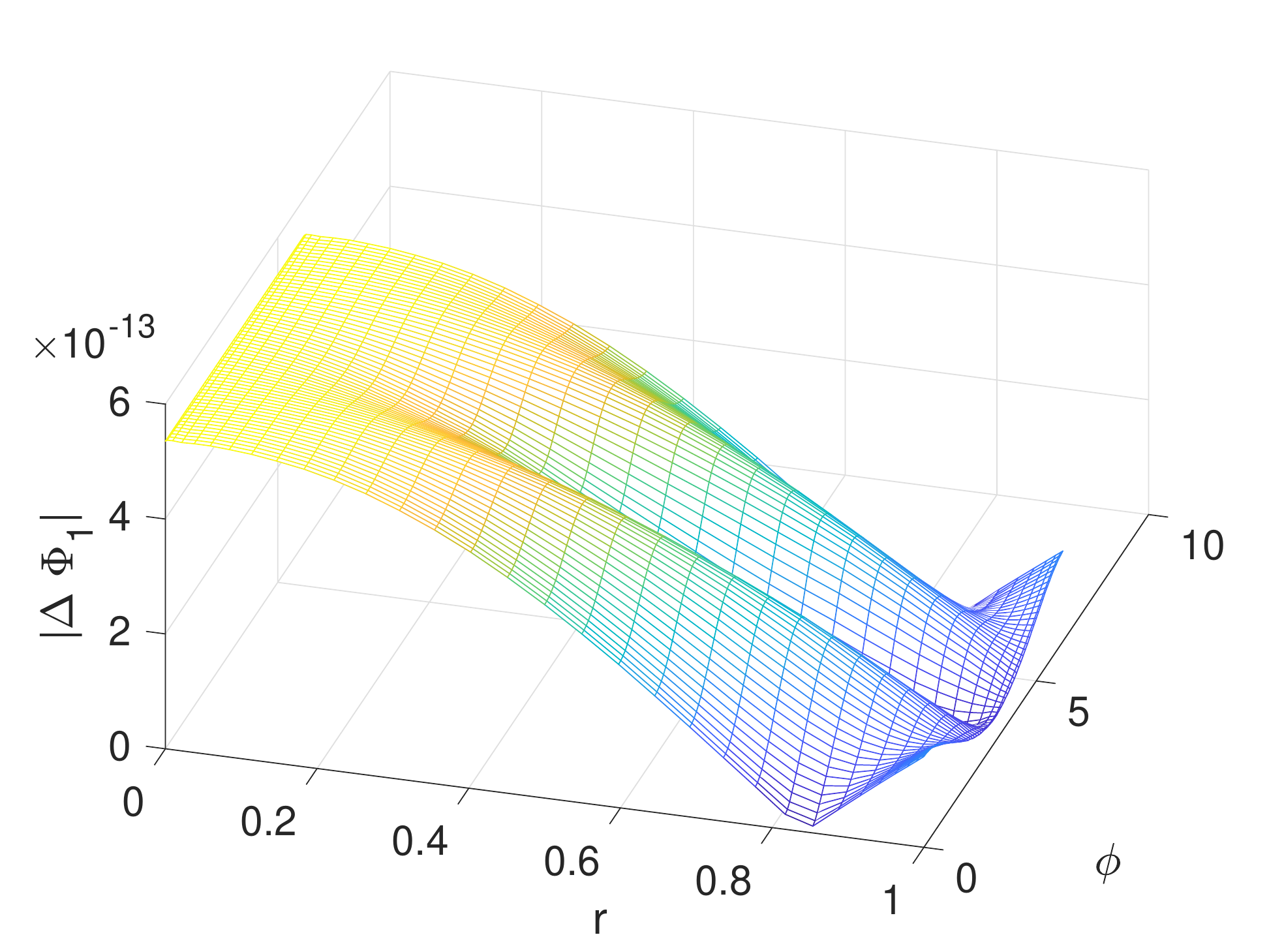}
   \includegraphics[width=0.49\textwidth]{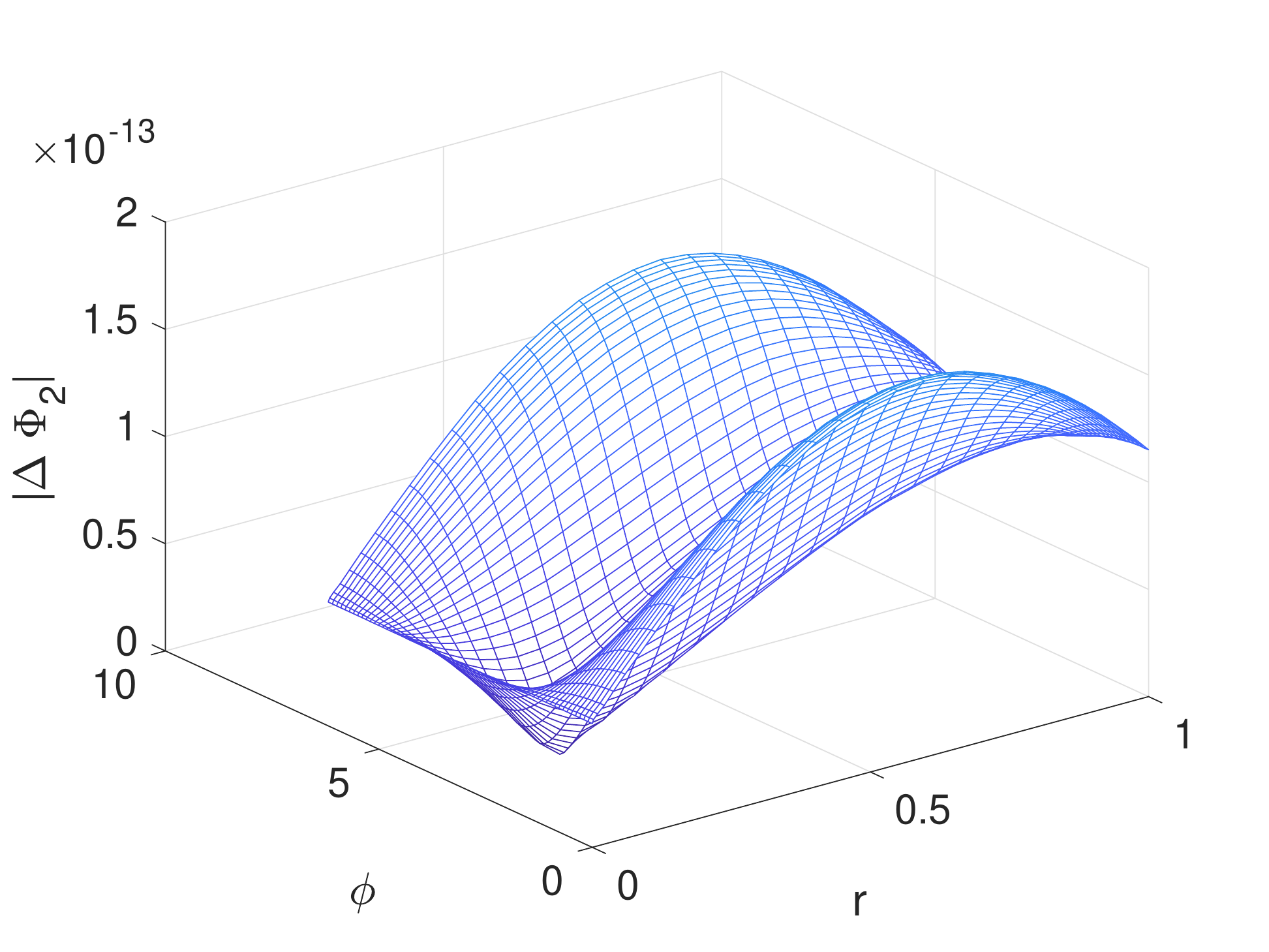}
 \caption{Differences of the CGO solutions to (\ref{dbarpsi}) and (\ref{dbarpsiasym}) 
 for $q=1$ on the unit disk and vanishing elsewhere for $k=1$, 
 obtained by the codes of the previous and the current section 
 respectively, on the 
 left for $\Phi_{1}$, on the right for $\Phi_{2}$.}
 \label{DeltaPhi}
\end{figure}

As in the previous sections, we consider the case of $q=1$ at the 
disk and vanishing outside. To test the iterative code, we compare 
its results for $k=1$ with the code of the previous section (recall 
that the fundamental solution itself was tested in section 3 by 
comparison with modified Bessel functions). We use $N_{r}=32$ and 
$N_{\phi}=64$ in both cases. The solutions $\psi_{1,2}$ are computed 
via the fundamental solution, the solutions $\Phi_{1,2}$ with the 
iterative code. The differences between the solutions 
$\Delta\Phi_{1}:=\mathrm{e}^{kz}\psi_{1}-\Phi_{1}$ and 
$\Delta\Phi_{2}:=\mathrm{e}^{\bar{k}\bar{z}}\psi_{2}-\Phi_{2}$ can be 
seen in Fig.~\ref{DeltaPhi}. They are both of the order $10^{-13}$ 
which implies that both solutions are determined with the same 
precision. 

\begin{figure}[htb!]
  \includegraphics[width=0.49\textwidth]{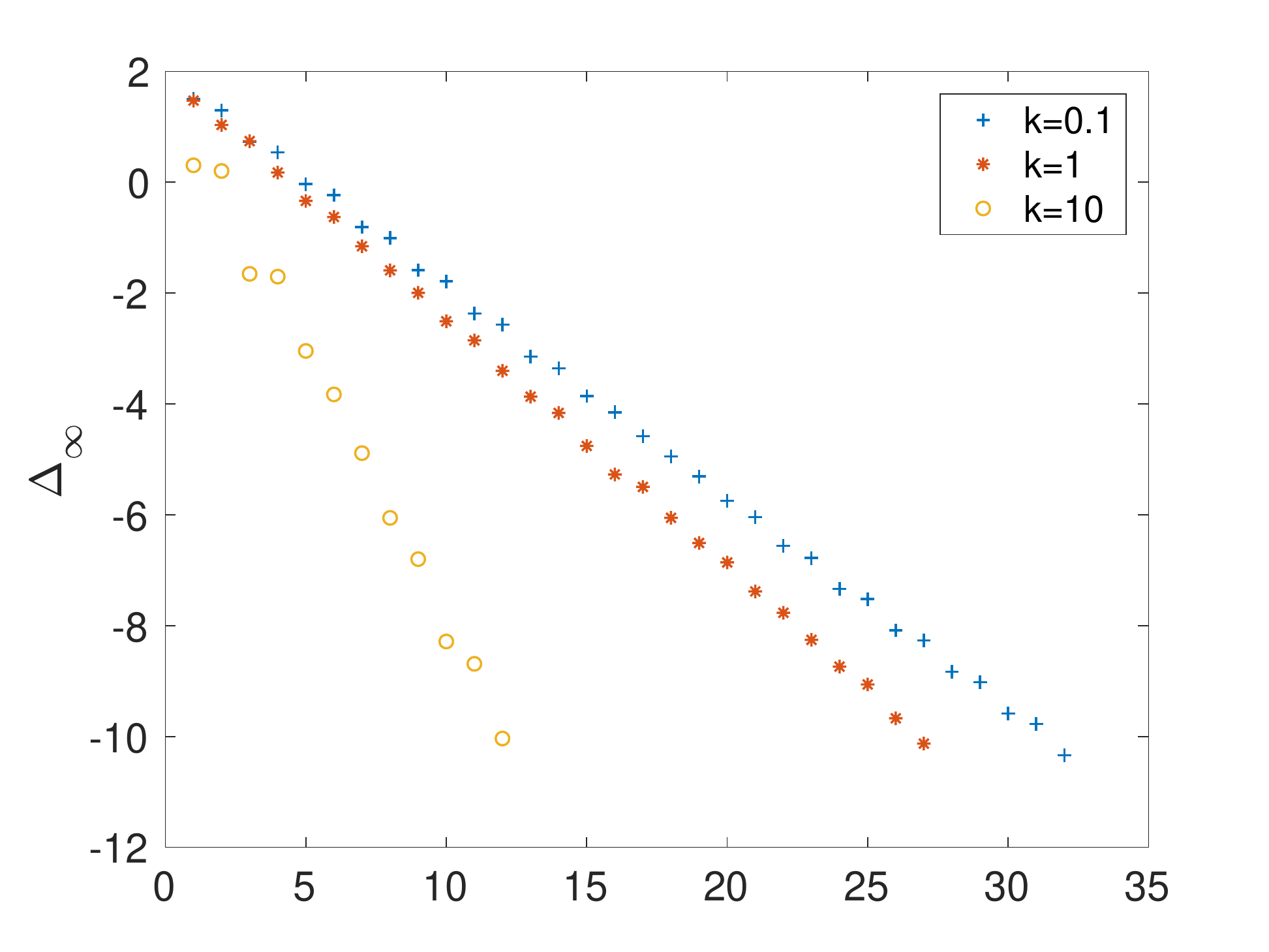}
  \includegraphics[width=0.49\textwidth]{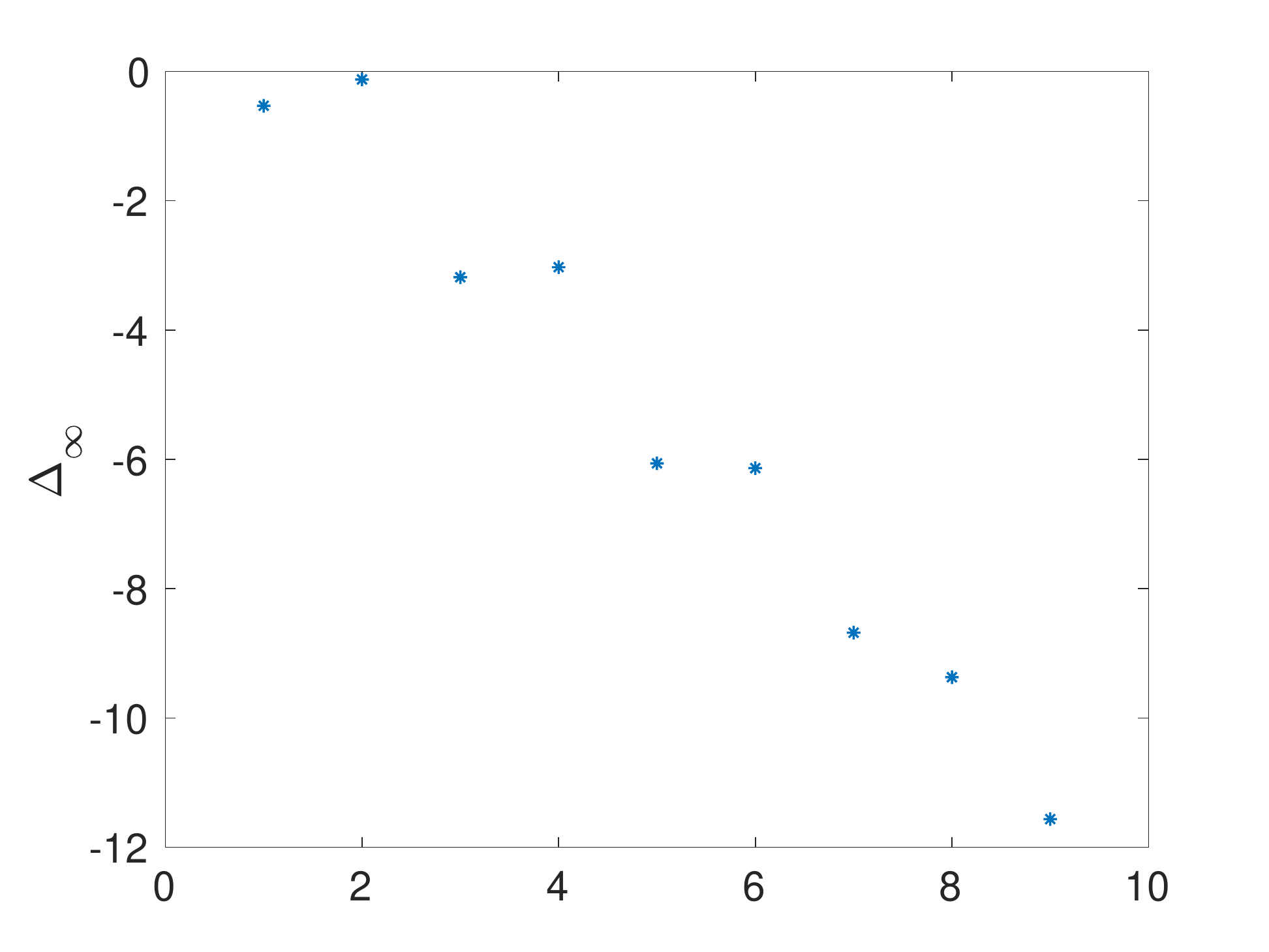}
 \caption{The quantity 
 $\Delta_{\infty}:=||\Phi_{1}^{(j+1)}-\Phi_{1}^{(j)}||_{\infty}+||\Phi_{2}^{(j+1)}-\Phi_{2}^{(j)}||_{\infty}$ 
in dependence of the number of iterations $j$ for $q=1$ on the unit 
disk and vanishing elsewhere, on the left  for $k=0.1,1,10$ and on 
the right for $k=100$.}
 \label{itconv}
\end{figure}

The iteration is generally stopped when 
$\Delta_{\infty}:=||\Phi_{1}^{(j+1)}-\Phi_{1}^{(j)}||_{\infty}+||\Phi_{2}^{(j+1)}-\Phi_{2}^{(j)}||_{\infty}$ 
is smaller than $10^{-10}$. For small $|k|$, this takes roughly 30 
iterations, but the iteration converges faster for larger $|k|$. We 
study this for three values of $k$, $k=0.1,1,10$ for the example 
$q=1$ on the disk. We use $N_{r}=32$ and $N_{\phi}=128$. The quantity 
$\Delta_{\infty}$ can be seen for these cases on the left of 
Fig.~\ref{itconv}. Visibly the convergence is linear. For $k=100$, 
just 9 iterations are needed as can be seen on the right of 
Fig.~\ref{itconv}. For the latter computation, $N_{r}=200$ and 
$N_{\phi}=600$ were used.  Note that the convergence depends also on 
the norm $||q||_{\infty}$ which is here equal to 1. If much larger 
values are to be considered as in \cite{AKMM,KMS}, it might be 
necessary to solve the system (\ref{dbarsysnumit}) without iteration, 
i.e., as in the previous section by inverting a large matrix. We do 
not address this possibility here since it was not necessary for the 
studied examples. 

\begin{figure}[htb!]
  \includegraphics[width=0.49\textwidth]{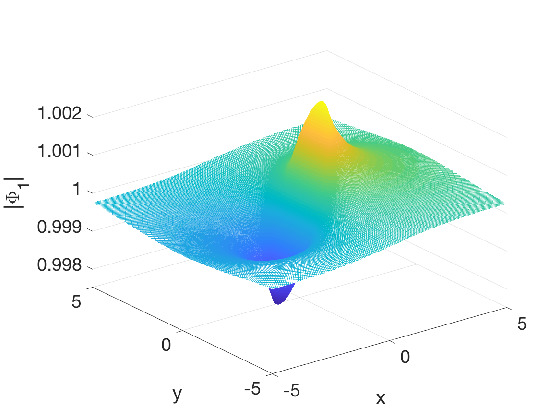}
  \includegraphics[width=0.49\textwidth]{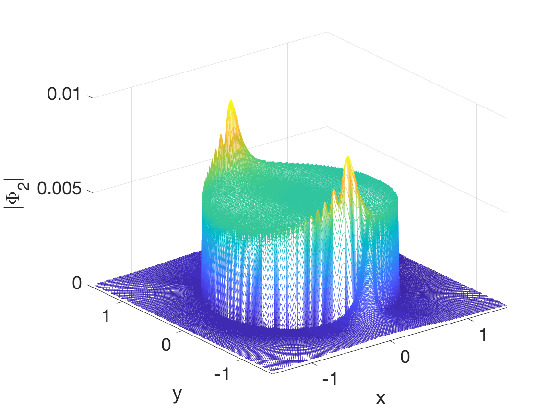}
 \caption{
 CGO solutions to (\ref{dbarpsi}) and (\ref{dbarpsiasym}) 
 for $q=1$ on the unit disk and vanishing elsewhere for $k=100$, on the 
 left for $\Phi_{1}$, on the right for $\Phi_{2}$.}
 \label{Phik100}
\end{figure}

The modulus of the solutions for $k=100$ can be seen in 
Fig.~\ref{Phik100}. The function $\Phi_{1}$ on the left of the figure 
appears to be roughly equal to 1 in correspondence with its 
asymptotic value with corrections of order $1/k$. The high frequency 
oscillations of the solution are hardly visible. These oscillations 
are much more pronounced for function $\Phi_{2}$ on the right of the 
same figure with amplitude of order $1/k$ around the asymptotic value 
0 for the solution. 

\begin{figure}[htb!]
  \includegraphics[width=0.49\textwidth]{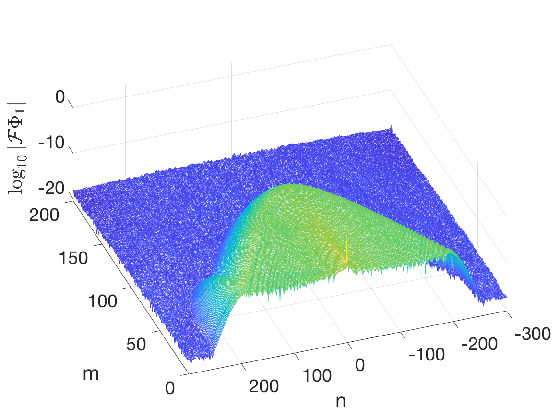}
  \includegraphics[width=0.49\textwidth]{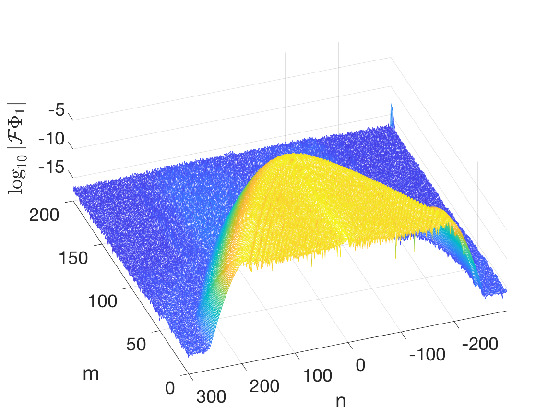}
 \caption{Spectral coefficients of CGO solutions to (\ref{dbarpsi}) and (\ref{dbarpsiasym}) 
 for $q=1$ on the unit disk and vanishing elsewhere for $k=100$, on the 
 left for $\Phi_{1}$, on the right for $\Phi_{2}$.}
 \label{Phik100coeff}
\end{figure}

As before, the numerical resolution can be checked via the decrease of 
the spectral coefficients for large value of the Chebyshev and 
Fourier indices. For the example $k=100$, this can be seen in 
Fig.~\ref{Phik100coeff}. The coefficients decrease as expected to the 
order of machine precision which shows that the solutions are 
numerically well resolved.

\begin{figure}[htb!]
  \includegraphics[width=0.32\textwidth]{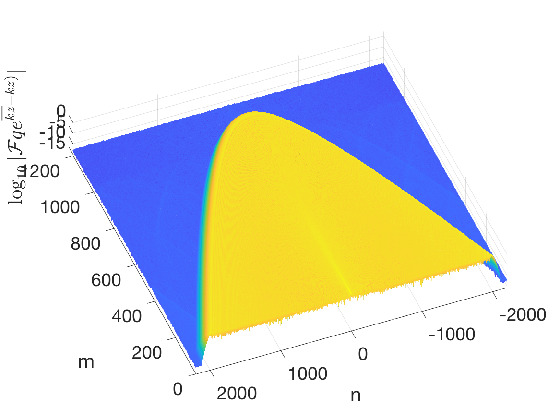}
  \includegraphics[width=0.32\textwidth]{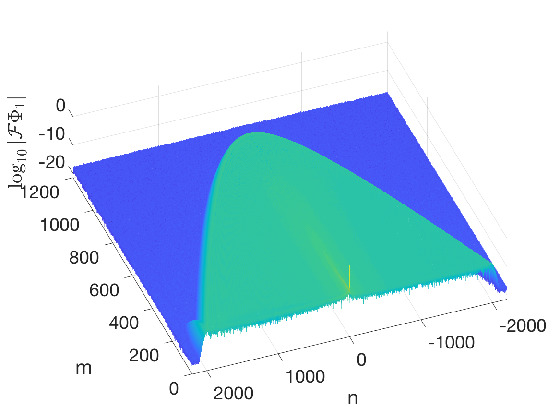}
  \includegraphics[width=0.32\textwidth]{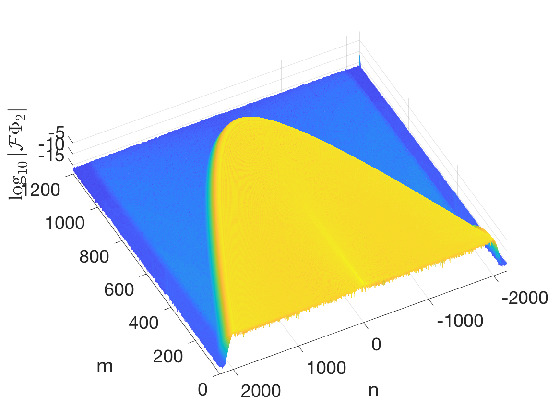}
 \caption{
 Spectral coefficients of $\mathrm{e}^{\bar{k}\bar{z}-kz}q$ 
 on the left, and  of CGO solutions to (\ref{dbarpsi}) and (\ref{dbarpsiasym}) 
 for $q=1$ on the unit disk and vanishing elsewhere for $k=1000$, on the 
 middle for $\Phi_{1}$, on the right for $\Phi_{2}$.}
 \label{qk1000}
\end{figure}

Note that no aliasing problems are observed here, which is due exactly to the 
fact that the spectral coefficients decrease both in Chebyshev and 
Fourier indices to machine precision. In fact much higher values of 
$|k|$ can be treated if this is respected. As discussed in remark 
\ref{remresolution}, the number $N_{r}$ of Chebyshev polynomials and 
the number $N_{\phi}$ of Fourier modes have to be chosen such that 
the spectral coefficients of $q(r,\phi)\mathrm{e}^{\bar{k}\bar{z}-kz}$ 
decrease to machine precision. For $k=1000$ this is shown in 
Fig.~\ref{qk1000} on the left. For $N_{r}=1200$ and 
$N_{\phi}=4400$, the spectral coefficients decrease to the order of 
the rounding error. With this choice of the parameter $N_{r}$ and $N_{\phi}$, the 
iteration converges after just 7 steps. The spectral 
coefficients shown in Fig.~\ref{qk1000} decrease as expected 
to the order of the rounding error. 

\begin{figure}[htb!]
  \includegraphics[width=0.49\textwidth]{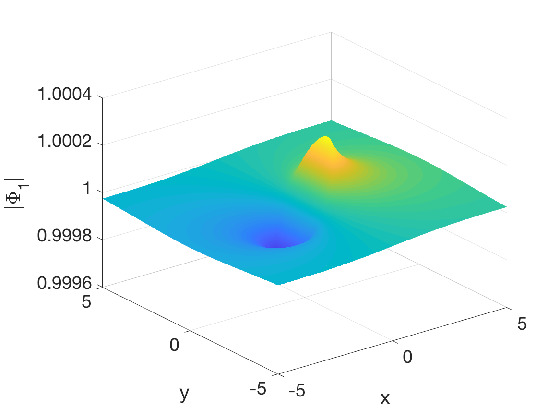}
  \includegraphics[width=0.49\textwidth]{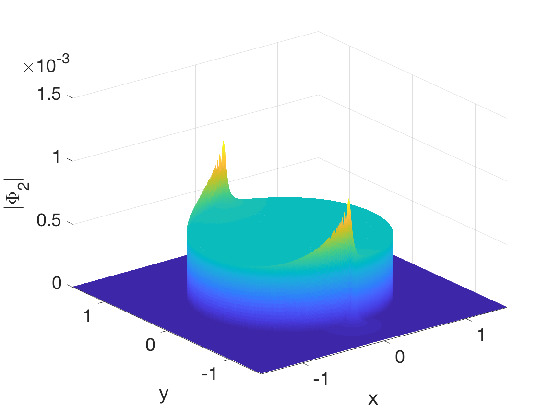}
 \caption{
 CGO solutions to (\ref{dbarpsi}) and (\ref{dbarpsiasym}) 
 for $q=1$ on the unit disk and vanishing elsewhere for $k=1000$, on the 
 left for $\Phi_{1}$, on the right for $\Phi_{2}$.}
 \label{Phi1k1000}
\end{figure}

The solutions can be 
seen in Fig.~\ref{Phi1k1000}. The function $\Phi_{1}$ on the left has 
hardly visible oscillations, the deviation from the asymptotic value 
1 is of order $1/k$ as in Fig.~\ref{Phik100}. The function $\Phi_{2}$ 
on the right on the other hand shows rapid oscillations of order 
$1/k$ around 0. 

\subsection{Reflection coefficient}

\begin{figure}[htb!]
  \includegraphics[width=0.49\textwidth]{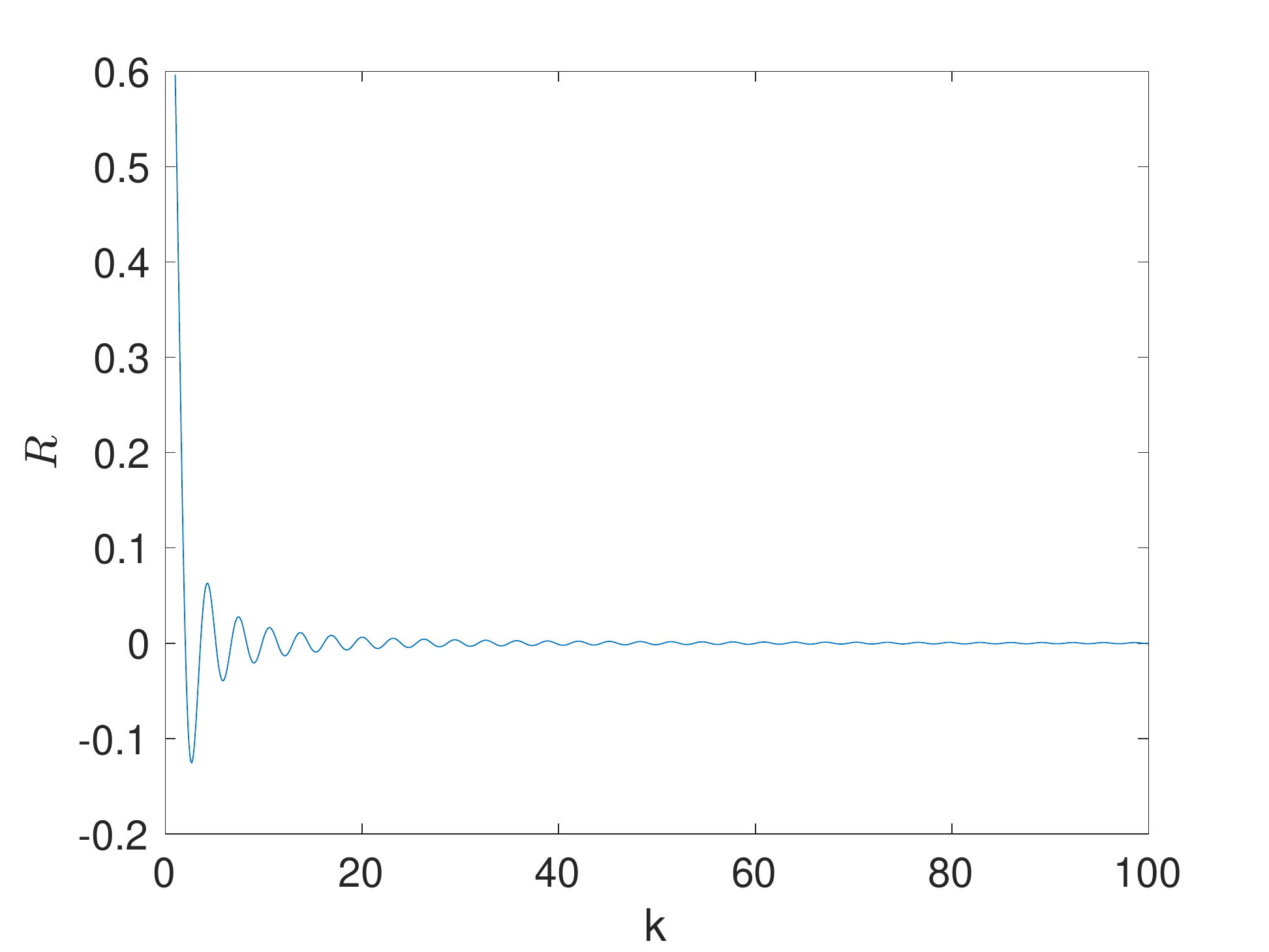}
  \includegraphics[width=0.49\textwidth]{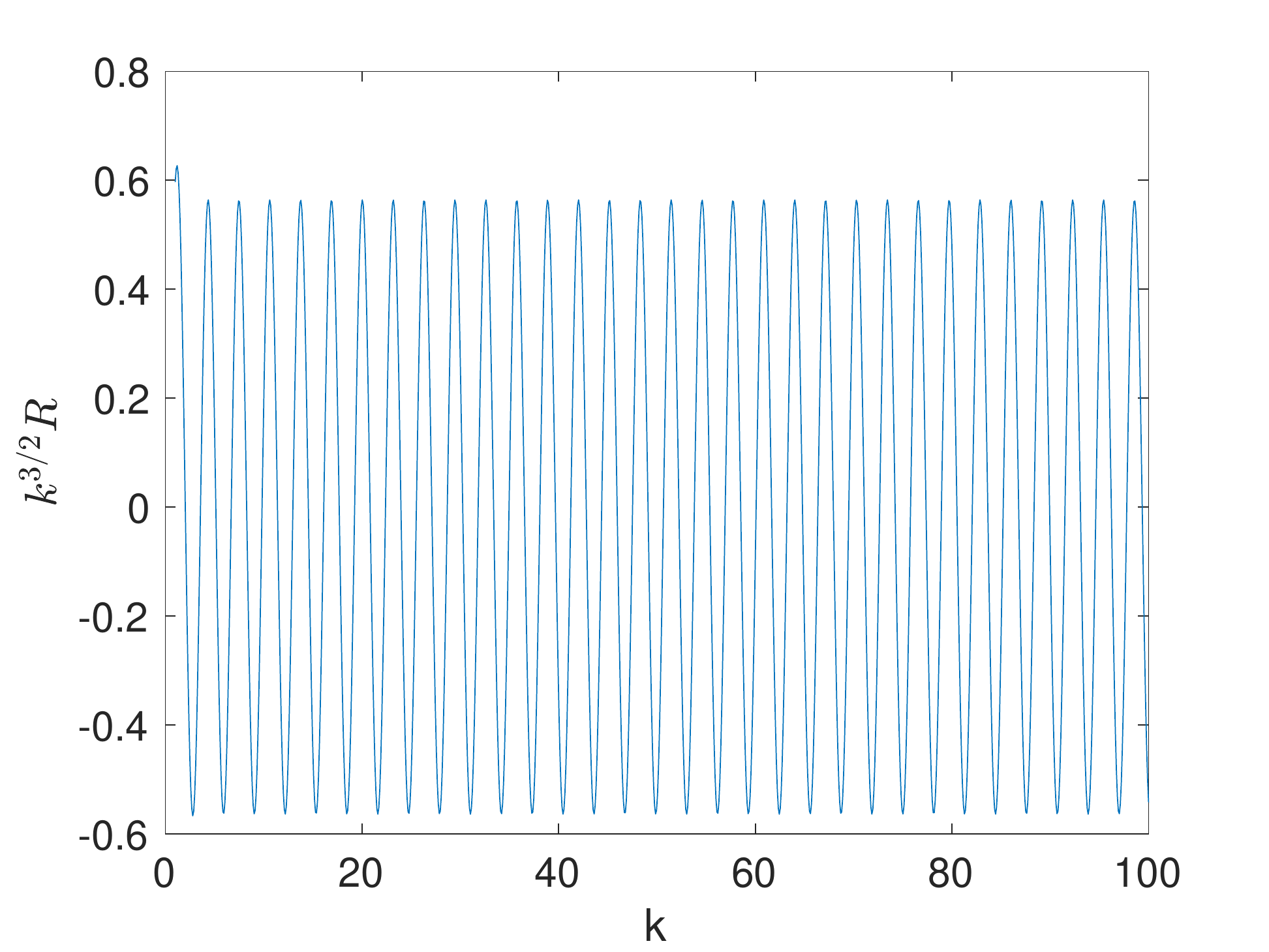}
 \caption{Reflection coefficient $R$ (\ref{eq:r-def})  
 for $q=1$ on the unit disk and vanishing elsewhere in dependence of 
 $k$ on the 
 left, and $kR$ on the right.}
 \label{reflectionk100}
\end{figure}

The reflection coefficient can be computed for a given value of $k$ 
via (\ref{reflection}). For the example of the characteristic 
function of the disk studied here, one gets for $k\in[1,100]$ the 
left figure of Fig.~\ref{reflectionk100}. The reflection coefficient has  an 
amplitude decreasing proportional $1/|k|^{3/2}$ and an oscillatory singularity at 
infinity as  can be seen from 
the plot of $Rk^{3/2}$ on the right of the same figure. 

An asymptotic formula for the reflection coefficent can be obtained 
via (\ref{eq:r-def}) and (\ref{eq:solid-Cauchy}) with 
$k=\kappa e^{i\psi}$, $\kappa,\psi\in\mathbb{R}$,
\begin{equation}
    \bar{R}(k)=\frac{2}{\pi}\int_{0}^{1}\int_{0}^{2\pi}q(r,\phi)e^{2i\kappa 
    r\sin(\phi-\psi)}\Phi_{1}rdr\phi.
    \label{Rasym}
\end{equation}
In \cite{joh} it was shown for $C^{\infty}$ potentials such that all 
derivatives are in $L^{2}(\mathbb{R}^{2})$ that $\Phi_{1}$ and 
$\Phi_{2}$ converge like a geometric series in $1/\kappa$ for $\kappa$ 
large. The situation is less clear for discontinuous potentials as 
the ones studied here, but numerical results indicate that 
$\Phi_{1}=1+O(1/\kappa)$ also in this case.  This would imply for large 
$\kappa$ with a stationary phase 
approximation for (\ref{Rasym}) for the example of the disk (where 
$k$ can be chosen to be real) that
\begin{equation}
    \bar{R}(k) \approx \frac{1}{\sqrt{\pi k^{3}}}\cos(2k-3\pi/4)
    \label{Rasym1}.
\end{equation}

\begin{figure}[htb!]
  \includegraphics[width=0.49\textwidth]{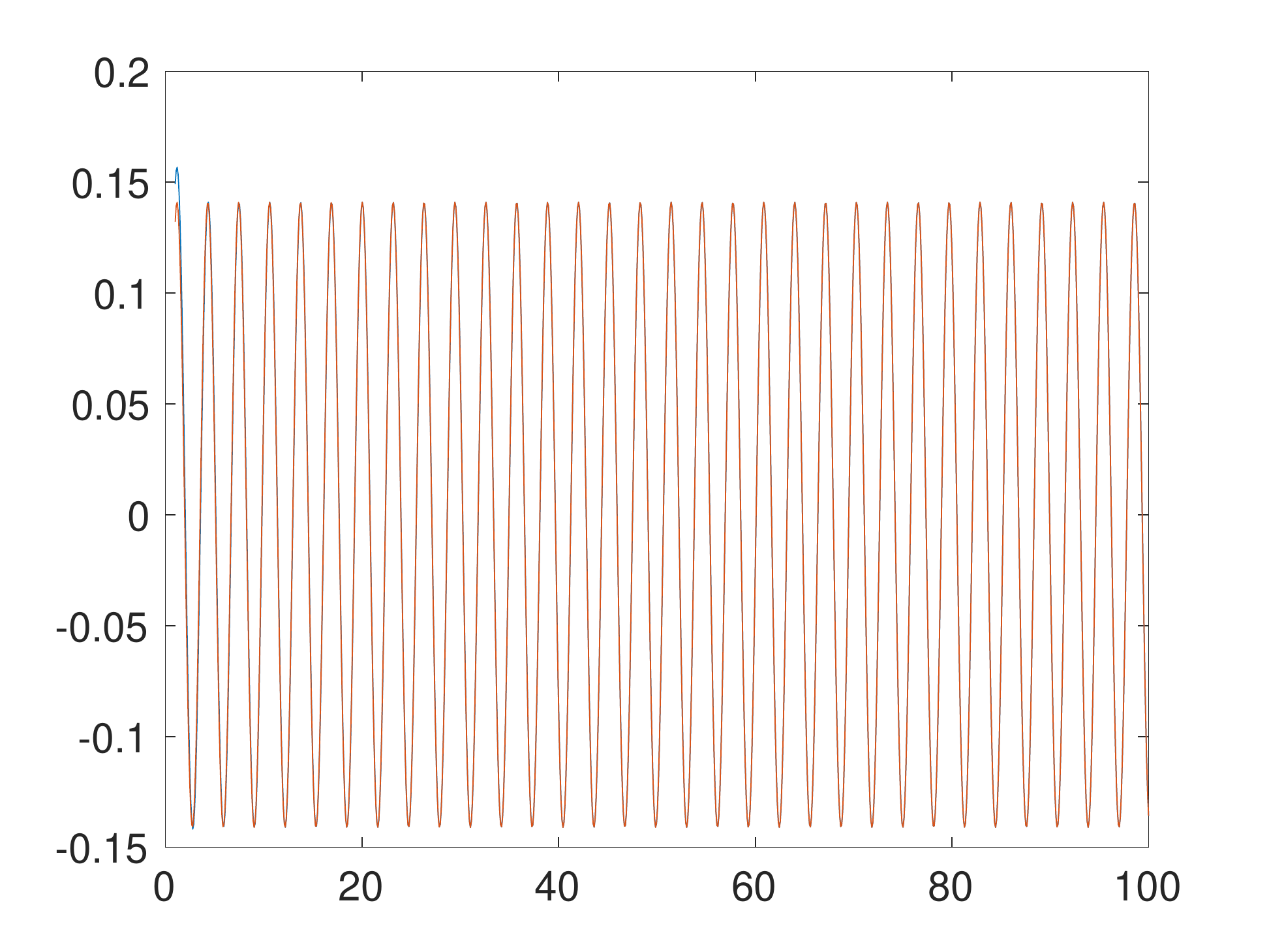}
  \includegraphics[width=0.49\textwidth]{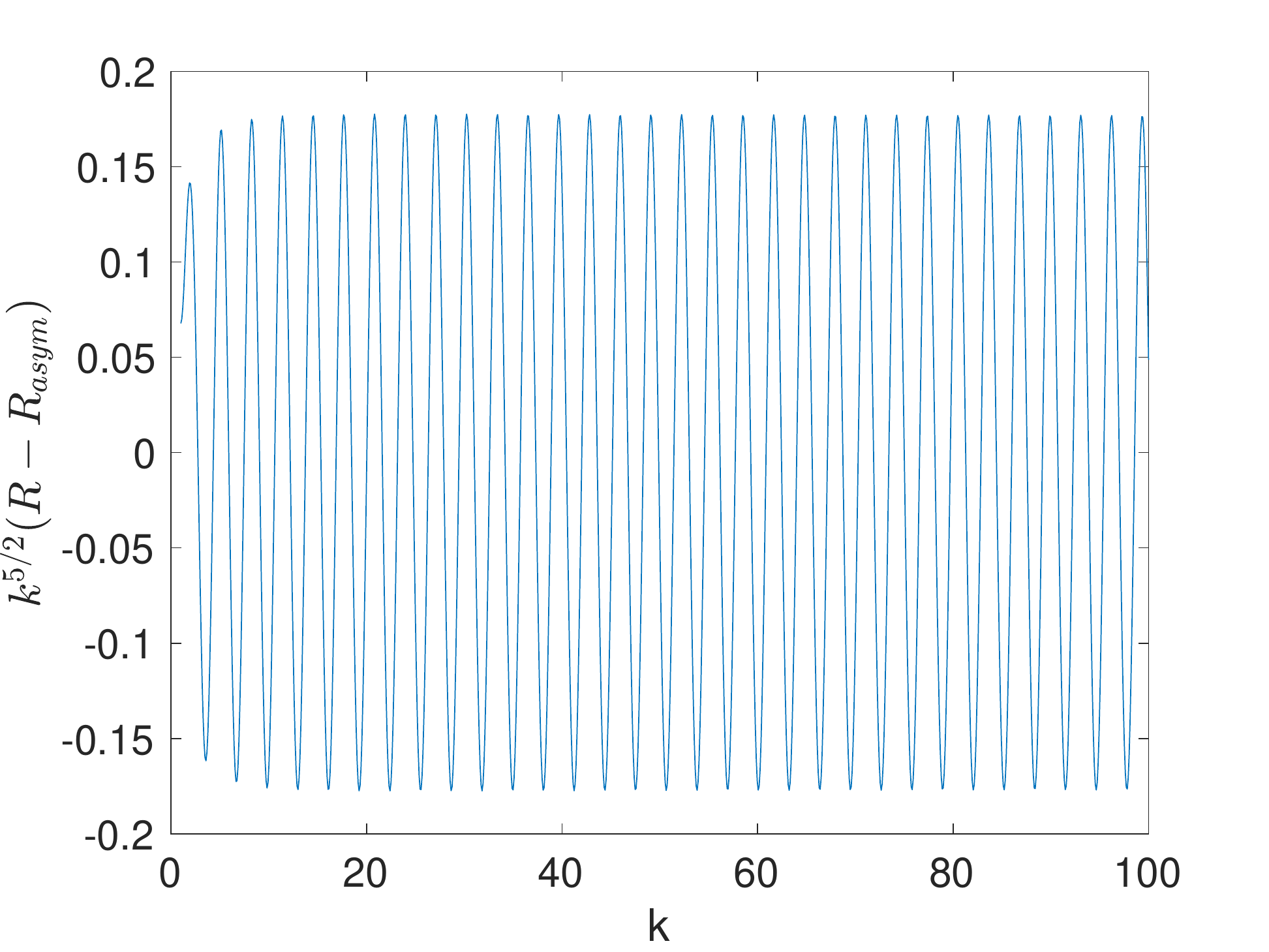}
 \caption{Reflection coefficient $R$ (\ref{eq:r-def})  
 for $q=1$ on the unit disk and vanishing elsewhere in dependence of 
 $k$ on the 
 left, in blue the numerical value, in red the asymptotic formula 
 $R_{asym}$
 (\ref{Rasym1}), and $k^{5/2}(R-R_{asym})$ on the right.}
 \label{reflectionkasym}
\end{figure}

In Fig.~\ref{reflectionkasym} we show on the left the numerically 
computed reflection coefficient in blue and in red the asymptotic 
formula (\ref{Rasym1}). It can be seen that the agreement is 
excellent already for values of the order of $|k|\sim 10$. The 
difference between numerically computed reflection coefficient and 
the asymptotic formula is shown to be of the order of $k^{5/2}$ on 
the right. 

\section{Outlook}
In this paper, we have presented spectral approaches for the numerical 
construction of CGO solutions for potentials with compact support on 
a disk, i.e.,  the scattering for the defocusing DS II equation. It is straight forward to generalize these approaches to 
functions which are piecewise smooth on annular parts of the disk. 
For large values of the modulus of the spectral parameter $|k|$, the 
reflection coefficient is shown to have an oscillatory behavior. An 
asymptotic formula is given in this case. This could allow the computation 
of the reflection coefficient, i.e., the scattering data for all 
values of $k$ with a \emph{hybrid} approach: for values of 
$|k|\lesssim 1000$, the reflection coefficient is numerically computed 
to machine precision. For larger values of $|k|$, an asymptotic series 
in $1\sqrt{|k|}$ will be given which could also  be computed for large $|k|$  
to machine precision. 
The oscillatory nature of the reflection coefficient for $k\to\infty$ 
makes the use of hybrid approaches necessary in the numerical 
inverse scattering. This means 
that the main asymptotic contribution has to be treated analytically with 
some semiclassical techniques, and just the residual between 
asymptotic and full CGO solution will be computed numerically. This 
will be the subject of further research. 

A further direction of research will be related to the exceptional 
points of the d-bar system (\ref{dbarpsi}) in the focusing case 
($\sigma=-1$). The numerical tools developed in the present paper and in 
\cite{KMS} will allow to study large classes of potentials to see in 
which form and when exceptional points appear, and to which physical 
features (appearence of lump solitons?) they correspond in the 
context of the focusing DS. 

An interesting question in the context of potentials with compact 
support would be cases with a non-circular boundary. The task would 
be to map for instance a smooth and convex boundary with a conformal 
transformation to a circle. The construction of such tansformations, 
see for instance \cite{HPT}, and its combination with the codes 
presented here will be a direction of future research. 

The present paper presents a purely spectral approach to the d-bar 
system (\ref{dbarpsi}), i.e., all functions are expanded in terms of 
combined Chebyshev and Fourier series which are then truncated. 
Actual computations are done with the coefficients of these series 
called \emph{spectral coefficients}. Resolution in terms of Chebychev 
polynomials is limited by the conditioning (second order 
differentiation matrices are 
known to have a conditioning of $O(N_{r}^{4})$, see the discussion in 
\cite{trefethen}). In \cite{OT} an ultraspherical approach with 
better conditioning and sparser matrices was presented which was applied in the context of 
the hypergeometric equation in \cite{CFKSV} to the 2d Laplace equation in polar 
coordinates. It was shown there that higher accuracies can be 
achieved than with  Chebychev differentiation matrices. We will test 
in future publications whether this approach can also be efficient in 
the present context. 

Note that the presented codes are highly parallelisable. First the task to 
compute CGO solutions for a high number of values of the spectral 
parameter $k$ can be done in parallel, in an efficient way on low 
cost GPUs as in \cite{KS}. Since the code is spectral and thus very 
efficient, they can be run on a single GPU for each $k$ up to rather 
large values of $|k|$. For exceptionally high resolutions, the codes 
can also be parallelised for single $k$ computation.

\section*{Acknowledgement}
This work was partially supported by the PARI and FEDER programs in 
2016 and 2017, by the ANR-FWF project ANuI, the isite BFC project 
NAANoD and by the Marie-Curie 
RISE network IPaDEGAN. We are indebted to J. Sj\"ostrand for 
countless helpful discussions and hints.

\end{document}